\documentclass{llncs}

\usepackage{times,amsmath,amssymb,enumerate,gastex}

\newcommand{\rest}{\mathord\restriction}
\newcommand{\NF}{\mathrm{NF}}
\newcommand{\RAT}{\mathrm{RAT}}
\newcommand{\bigast}{\mathop{\hbox{\Large $\ast$}}}

\begin{document}

\title{The submonoid and rational subset membership problems for graph groups}
\author{Markus Lohrey\inst{1} \and Benjamin Steinberg\inst{2,}\thanks{The second author would like to
acknowledge the support of an NSERC grant.} \institute{Universit\"at
Stuttgart, FMI, Germany \and School of Mathematics and Statistics,
Carleton University, ON, Canada\\
\email{lohrey@informatik.uni-stuttgart.de,
bsteinbg@math.carleton.ca}}}

\maketitle

\begin{abstract}
We show that the membership problem in a finitely generated
submonoid of a graph group (also called a right-angled Artin group
or a free partially commutative group) is decidable if and only if
the independence graph (commutation graph) is a transitive forest.
As a consequence we obtain the first example of a finitely presented
group with a decidable generalized word problem that does not have a
decidable membership problem for finitely generated submonoids.  We
also show that the rational subset membership problem is decidable
for a graph group if and only if the independence graph is a
transitive forest, answering a question of Kambites, Silva, and the
second author \cite{KaSiSt06}.  Finally we prove that for certain 
amalgamated free products and HNN-extensions the rational subset and
submonoid membership problems are recursively equivalent.
In particular, this applies to finitely generated groups with two or
more ends that are either torsion-free or residually finite.
\end{abstract}

\section{Introduction}

Algorithmic problems concerning groups are a classical topic in
algebra and theoretical computer science. Since the pioneering work
of Dehn from 1910 \cite{Dehn10}, decision problems like the word
problem or the generalized word problem
(which is also known as the subgroup membership problem since it asks
whether one can decide if a given group element belongs to a given
finitely generated subgroup) have been intensively studied for
various classes of groups. A first natural generalization of these
classical decision problems is the submonoid membership problem:
given a finite set $S$ of elements of $G$ and an element $g\in G$,
does $g$ belong to the submonoid generated by $S$?
Notice that $g$ has finite order if and only if $g^{-1}$ is in the
submonoid generated by $g$ and so decidability of the submonoid
membership problem lets one determine algorithmically the order of
an element of the group $G$.  A recent paper on the submonoid
membership problem is Margolis, Meakin, and \v Suni\'k~\cite{MaMeSu}.

A further generalization is the rational subset membership problem:
for a given rational subset $L$ of a group $G$ and an element $g \in
G$ it is asked whether $g \in L$. The class of rational subsets of a
group $G$ is the smallest class that contains all finite subsets of
$G$, and which is closed under union, product, and the Kleene hull
(or Kleene star; it associates to a subset $L \subseteq G$ the submonoid $L^*$
generated by $L$).  Equivalently, it consists of the all subsets of
$G$ recognizable by finite automata.   Rational subsets in arbitrary
groups and monoids are an important research topic in language theory,
see, e.g., \cite{Ber79,KaSiSt06,Ned00}. The rational subset membership
problem generalizes the 
submonoid membership problem and the the generalized word problem
for a group, because every finitely generated submonoid (and hence
subgroup) of a group is rational.

It is easy to see that decidability of the rational subset
membership problem transfers to finitely generated subgroups.
Grunschlag has shown that the property of having a decidable
rational subset membership problem is preserved under finite
extensions, i.e., if $G$ has a decidable rational subset membership
problem and $G \leq H$, where the index of $G$ in $H$ is finite,
then  $H$ also has a decidable rational subset membership problem
\cite{Gru92}. Kambites, Silva, and the second author \cite{KaSiSt06}
proved that the fundamental group of a finite graph of
groups~\cite{Se03} with finite edge groups has a decidable rational
subset membership problem provided all vertex groups have a
decidable rational subset membership problem. In particular, this
implies that decidability of the rational subset membership problem
is preserved by free products, see also \cite{Ned00}.

The main result of this paper is to characterize the decidability of
the submonoid membership problem and the rational subset membership
problem for graph groups.  In particular we provide the first
example, as far as we know, of a group with a decidable generalized
word problem that does not have a decidable submonoid (and hence
rational subset) membership problem.

A \emph{graph group}~\cite{Dro85} $\mathbb{G}(\Sigma,I)$ is
specified by a finite undirected graph $(\Sigma,I)$, which is also
called an \emph{independence alphabet} (or \emph{commutation
graph}). The graph group $\mathbb{G}(\Sigma,I)$ is formally defined
as the quotient group of the free group generated by $\Sigma$ modulo
the set of all relations $ab=ba$, where $(a,b) \in I$. Graph groups
are a group analogue to trace monoids (free partially commutative
monoids), which play a prominent role in concurrency theory
\cite{DieRoz95}. Graph groups are also called \emph{free partially
commutative groups} \cite{Die90tcs,Wra88}, \emph{right-angled Artin
groups} \cite{BrMe01,CrWi04}, and \emph{semifree groups}
\cite{Bau81}. They are currently a hot topic of interest in group
theory, in particular because of the richness of the class of groups
embeddable in graph groups.  For instance, the Bestvina-Brady
groups, which were used to distinguish the finiteness properties
$\mathcal F_n$ and $\text{FP}_n$~\cite{BesBr97} (and were also
essential for distinguishing the finiteness properties FDT and FHT
for string rewriting systems~\cite{OtPr04}), are subgroups of graph
groups. Crisp and Wiest show that the fundamental group of any
orientable surface (and of most non-orientable surfaces) embeds in a
graph group~\cite{CrWi04}. Another class of groups that embed into
graph groups are fundamental groups of finite state complexes \cite{GhPe07}.

Algorithmic problems concerning graph groups have been intensively
studied in the past, see, e.g.,
\cite{Die90tcs,DiMu06,EsKaRe06,KaSiSt06,KaWeMy05,Wra88}. In
\cite{Die90tcs,Wra88} it was shown that the word problem for a graph
group can be decided in linear time (on a random access machine). A
recent result of Kapovich, Weidmann, and Myasnikov \cite{KaWeMy05}
shows that if $(\Sigma,I)$ is a chordal graph (i.e., if $(\Sigma,I)$
does not have an induced cycle of length at least 4), then the
generalized word problem for $\mathbb{G}(\Sigma,I)$ is decidable. On
the other hand, a classical result of Mihailova \cite{Mih66} states
that already the generalized word problem for the direct product of
two free groups of rank 2 is undecidable. Note that this group is
the graph group $\mathbb{G}(\Sigma,I)$, where the graph $(\Sigma,I)$
is a cycle on 4 nodes (also called $\mathsf{C4}$). In fact,
Mihailova proves a stronger result: she constructs a \emph{fixed}
subgroup $H$ of $\mathbb{G}(\mathsf{C4})$ such that it is
undecidable, whether a given element of $\mathbb{G}(\mathsf{C4})$
belongs to $H$. Recently, it was shown by Kambites that a graph
group $\mathbb{G}(\Sigma,I)$ contains a direct product of two free
groups of rank 2 if and only if $(\Sigma,I)$ contains an induced
$\mathsf{C4}$~\cite{Kam06}. This leaves a gap between the
decidability result of \cite{KaWeMy05} and the undecidability result
of Mihailova \cite{Mih66}.

In \cite{KaSiSt06} it is shown that the rational subset membership
problem is decidable for a free product of direct products of a free
group with a free Abelian group. Such a group is a graph group
$\mathbb{G}(\Sigma,I)$, where every connected component of
$(\Sigma,I)$ results from connecting all nodes of a clique with all
nodes from an edge-free graph. On the other hand, the only
undecidability result for the rational subset membership problem for
graph groups that was known so far is Mihailova's result for
independence alphabets containing an induced $\mathsf{C4}$.

In this paper, we shall characterize those graph groups for which
the rational subset membership problem is decidable: we prove that
these are exactly those graph groups $\mathbb{G}(\Sigma,I)$, where
$(\Sigma,I)$ is a transitive forest (Theorem~\ref{T main}). The
graph $(\Sigma,I)$ is a transitive forest if it is the disjoint
union of comparability  graphs of rooted trees. An alternative
characterization of transitive forests was presented in
 \cite{Wolk65}: $(\Sigma,I)$ is a transitive forest if and
 only if it neither contains an induced $\mathsf{C4}$
 nor an induced path on 4 nodes (also called $\mathsf{P4}$).
Graph groups $\mathbb{G}(\Sigma,I)$, where $(\Sigma,I)$ is a
transitive forest, have also appeared in \cite{MeRa06}: they are
exactly those graph groups which are subgroup separable (the case of
$\mathsf{P4}$ appears in~\cite{NiWi01}). Recall that a group $G$ is
called subgroup separable if, for every finitely generated subgroup
$H \leq G$ and every $g \in G \setminus H$ there exists a normal
subgroup $N \leq G$ having finite index such that 
$g \not\in NH$.  Subgroup separability implies decidability of the
generalized word problem.

One half of Theorem~\ref{T main} can be easily obtained from a
result of Aalbersberg and Hoogeboom \cite{AaHo89}: The problem of
deciding whether the intersection of two rational subsets of the
trace monoid (free partially commutative monoid)
$\mathbb{M}(\Sigma,I)$ is nonempty is decidable if and only if
$(\Sigma,I)$ is a transitive forest. Now, $L \cap K \neq \emptyset$
for two given rational subsets $L, K \subseteq \mathbb{M}(\Sigma,I)$
if and only if $1 \in L K^{-1}$ in the graph group
$\mathbb{G}(\Sigma,I)$. Hence, if $(\Sigma,I)$ is not a transitive
forest, then the rational subset membership problem for
$\mathbb{G}(\Sigma,I)$ is undecidable. In fact, we construct a fixed
rational subset $L \subseteq \mathbb{G}(\Sigma,I)$ such that it is
undecidable whether $g \in L$ for a given group element $g \in
\mathbb{G}(\Sigma,I)$.

The converse direction in Theorem~\ref{T main}
is an immediate corollary of our
Theorem~\ref{T dec}, which is one of the main group theoretic
results of this paper.  It states that the rational subset
membership problem is decidable for every group that can be built up
from the trivial group using the following four operations: (i)
taking finitely generated subgroups, (ii) finite extensions, (iii)
direct products with $\mathbb{Z}$, and (iv) finite graphs of groups
with finite edge groups. Note that the only operation that is not
covered by the results cited earlier is the direct product with
$\mathbb{Z}$. In fact, it seems to be an open question whether
decidability of the rational subset membership problem is preserved
under direct products with $\mathbb{Z}$. Hence, we have to follow
another strategy. We will introduce a property of groups that
implies the decidability of the rational subset membership problem,
and which has all the desired closure properties. Our proof of
Theorem~\ref{T dec} uses mainly techniques from formal language theory
(e.g., semilinear sets, Parikh's theorem) and is inspired by
the methods from \cite{AaHo89,Bou94}.

It should be noted that due to the above reduction from the
intersection problem for rational trace languages to the rational
subset membership problem for the corresponding graph group, we also
obtain an alternative to the quite difficult proof from
\cite{AaHo89} for the implication ``$(\Sigma,I)$ is a transitive
forest $\Rightarrow$ intersection problem for rational subsets of
$\mathbb{M}(\Sigma,I)$ is decidable''.

In Section~\ref{S submonoid} we consider the \emph{submonoid
membership problem} for groups. We prove that for an amalgamated free
product $G \ast_A H$ such that $A$ is a finite proper subgroup of $G$
and $H$ and  
there exist $g \in G$, $h \in H$ with $g^{-1} A g \cap A = 1 = h^{-1} A h \cap A$,
the rational subset membership
problem is recursively equivalent to the  submonoid membership
problem (Theorem~\ref{T submonoid}). An analogous result is proved for
certain HNN extensions with finite associated subgroups.
%The idea is to encode the
%states of a finite automaton into one of the two factors of $G \ast_A H$.
As a consequence we obtain that the rational subset
membership problem is recursively equivalent to the  submonoid
membership problem for a group with two or more ends that is either
torsion-free or residually finite (Corollary~\ref{C tf or rs}).  
Using similar techniques, we are
also able to prove that the submonoid membership problem is
undecidable for the graph group $\mathbb{G}(\Sigma,I)$, where
$(\Sigma,I)$ is $\mathsf{P4}$ (Theorem~\ref{T P4}).   The result
of~\cite{KaWeMy05} shows that this graph group does have a decidable
generalized word problem, thereby giving our example of a group with
a decidable generalized word problem but an undecidable submonoid
membership problem. Together with Mihailova's undecidability result
for $\mathsf{C4}$ and our decidability result for transitive forests
(Theorem~\ref{T main}) it also follows that the submonoid membership
problem for a graph group $\mathbb{G}(\Sigma,I)$ is decidable if and
only if $(\Sigma,I)$ is a transitive forest (Corollary~\ref{C
submonoid graph group}).

Another consequence of our results is that the rational subset
membership problem for groups is recursively equivalent to the
submonoid membership problem if and only if a free product of groups
with decidable submonoid membership problems has a decidable submonoid
membership problem.

\section{Preliminaries}

We assume that the reader has some basic knowledge
in formal language theory (see, e.g., \cite{Ber79,HoUl79})
and group theory (see, e.g., \cite{LySch77,Rot95}).

\subsection{Formal languages}

Let $\Sigma$ be a finite alphabet. We use $\Sigma^{-1} = \{ a^{-1}
\mid a \in \Sigma\}$ to denote a disjoint copy of $\Sigma$. Let
$\Sigma^{\pm 1} = \Sigma \cup \Sigma^{-1}$. Define $(a^{-1})^{-1} =
a$; this defines an involution ${}^{-1} : \Sigma^{\pm 1} \to
\Sigma^{\pm 1}$, which can be extended to the free monoid $(\Sigma^{\pm 1})^*$ by
setting $(a_1 \cdots a_n)^{-1} = a_n^{-1} \cdots a_1^{-1}$. For a
word $w \in \Sigma^*$ and $a \in \Sigma$ we denote by $|w|_a$ the
number of occurrences of $a$ in $w$. For a subset $\Gamma \subseteq
\Sigma$, we denote by $\pi_{\Gamma}(w)$ the projection of the word
$w$ to the alphabet $\Gamma$, i.e., we erase in $w$ all symbols from
$\Sigma \setminus \Gamma$.

Let $\mathbb{N}^\Sigma$ be the set of all mappings from $\Sigma$ to
$\mathbb{N}$. By fixing an arbitrary linear order on the alphabet
$\Sigma$, we may identify a mapping $f \in \mathbb{N}^\Sigma$ with a
tuple from $\mathbb{N}^{|\Sigma|}$. For a word $w \in \Sigma^*$,
the \emph{Parikh image} $\Psi(w)$ is defined as the mapping $\Psi(w) :
\Sigma \to \mathbb{N}$ such that $[\Psi(w)](a) = |w|_a$ for all $a
\in \Sigma$. For a language $L \subseteq \Sigma^*$, the Parikh image
is $\Psi(L) = \{ \Psi(w) \mid w \in L \}$. For a set $K \subseteq
\mathbb{N}^\Sigma$ and $\Gamma \subseteq \Sigma$ let
$\overline{\pi}_\Gamma(K) = \{ f\rest_{\Gamma} \in
\mathbb{N}^{\Gamma} \mid f \in K \}$, where $f\rest_{\Gamma}$
denotes the restriction of $f$ to $\Gamma$.  We also need a notation
for the composition of erasing letters and taking the Parikh image. So,
for $L \subseteq \Sigma^*$ and $\Gamma \subseteq \Sigma$,
let $\Psi_\Gamma(L) = \overline{\pi}_\Gamma(\Psi(L)) (=
\Psi(\pi_\Gamma(L)))$; it may be viewed as a subset of
$\mathbb{N}^{|\Gamma|}$. A special case occurs when $\Gamma =
\emptyset$. Then either $\Psi_\emptyset(L) = \emptyset$ (if $L =
\emptyset$) or $\Psi_\emptyset(L)$ is the singleton set consisting
of the unique mapping from $\emptyset$ to $\mathbb{N}$.

A subset $K \subseteq \mathbb{N}^k$ is said to be \emph{linear} if there are
$x, x_1, \ldots, x_\ell \in \mathbb{N}^k$ such that $K = \{ x +
\lambda_1 x_1 +\cdots+\lambda_\ell x_\ell \mid \lambda_1, \ldots,
\lambda_\ell \in \mathbb{N}\}$, i.e.\ $K$ is a translate of a
finitely generated submonoid of $\mathbb{N}^k$.
A \emph{semilinear} set is a finite union of linear sets.

Let $\mathsf{G} = (N,\Gamma, S, P)$ be a context-free grammar, where $N$
is the set of nonterminals, $\Gamma$ is the terminal alphabet,
$S \in N$ is the start nonterminal, and $P \subseteq N \times (N \cup \Gamma)^*$
is the finite set of productions.
For $u, v \in (N \cup \Gamma)^*$ we write
$u \Rightarrow_{\mathsf{G}} v$ if $v$ can be derived from $u$ by
applying a production from $P$.
For $A \in N$, we define $L(\mathsf{G},A) = \{ w \in \Gamma^* \mid A
\stackrel{*}{\Rightarrow}_{\mathsf{G}} w \}$
and $L(\mathsf{G}) = L(\mathsf{G},S)$.
Parikh's theorem states that
the Parikh image of a context-free language is semilinear \cite{Par66}.

We will allow a more general form of productions
in context-free grammars, where the right-hand side of
a production is a regular language over the alphabet
$N \cup \Gamma$. Such a production $A \to L$ represents
the (possibly infinite) set of productions $\{ A \to s \mid s\in L\}$.
Clearly, such an extended context-free
grammar can be transformed effectively into an equivalent
context-free grammar with only finitely many productions.

Let $M$ be a monoid. The set $\RAT(M)$ of all
\emph{rational subsets} of $M$ is the smallest subset
of $2^M$, which contains all finite subsets of $M$, and
which is closed under  union, product, and
Kleene hull (the Kleene hull $L^*$
of a subset $L \subseteq M$ is the submonoid
of $M$ generated by $L$). By Kleene's theorem,
a subset $L \subseteq \Sigma^*$ is rational if and only if
$L$ can be recognized by a finite automaton.
If $M$ is generated by the finite set
$\Sigma$ and $h : \Sigma^* \to M$
is the corresponding canonical monoid homomorphism, then
$L \in \RAT(M)$ if and only if $L = h(K)$ for
some $K \in \RAT(\Sigma^*)$.
In this case, $L$ can be specified by a finite
automaton over the alphabet $\Sigma$.
The rational subsets of the free commutative monoid
$\mathbb{N}^k$ are exactly the semilinear subsets of
$\mathbb{N}^k$ \cite{EiSchu69}.

\subsection{Groups}

Let $G$ be a finitely generated group and let
$\Sigma$ be a finite group generating set for $G$.
Hence, $\Sigma^{\pm 1}$ is a finite monoid generating
set for $G$ and there exists a canonical monoid
homomorphism $h: (\Sigma^{\pm 1})^* \to G$.
The language
$$
\text{WP}_{\Sigma}(G) = h^{-1}(1)
$$
is called the \emph{word problem} of $G$ with respect to $\Sigma$,
i.e., $\text{WP}_{\Sigma}(G)$ consists of all words over the
alphabet $\Sigma^{\pm 1}$ which are equal to $1$ in the group $G$.
It is well known and easy to see that if $\Gamma$ is another finite
generating set for $G$, then $\text{WP}_{\Sigma}(G)$ is decidable if
and only if $\text{WP}_{\Gamma}(G)$ is decidable.

The \emph{submonoid membership problem} for $G$ is
the following decision problem:

\medskip

\noindent
INPUT: A finite set of words $\Delta \subseteq
(\Sigma^{\pm 1})^*$ and a word $w \in (\Sigma^{\pm 1})^*$.

\noindent
QUESTION: $h(w) \in h(\Delta^*)$?

\medskip

\noindent
Note that the subset $h(\Delta^*) \subseteq G$ is the
submonoid of $G$ generated by $h(\Delta) \subseteq G$.
If we replace in the submonoid membership problem the finitely
generated submonoid $h(\Delta^*)$ by the finitely generated subgroup
$h((\Delta \cup \Delta^{-1})^*)$, then we obtain the \emph{subgroup
membership problem}, which is also known as the \emph{generalized
word problem} for $G$. This term is justified, since the word
problem is a particular instance, namely with $\Delta = \emptyset$.
A generalization of the submonoid
membership problem for $G$ is the \emph{rational subset membership
problem}:

\medskip

\noindent
INPUT: A finite automaton $A$ over the alphabet
$\Sigma^{\pm 1}$ and a word $w \in
 (\Sigma^{\pm 1})^*$.

\noindent
QUESTION: $h(w) \in h(L(A))$?

\medskip

\noindent Note that $h(w) \in h(L(A))$ if and only if $1 \in
h(w^{-1}L(A))$. Since $w^{-1}L(A)$ is again a rational language, the
rational subset membership problem for $G$ is recursively equivalent
to the decision problem of asking whether $1 \in h(L(A))$ for a
given finite automaton $A$ over the alphabet $\Sigma^{\pm 1}$.
It should be noted that for all the computational problems
introduced above the decidability is
independent of the chosen generating set for $G$.

In the rational subset (resp. submonoid) membership problem, the
rational subset (resp. submonoid) is part of the input. Non-uniform
variants of these problems, where the
rational subset (resp. submonoid) is fixed, have been studied as well. More
generally, we can define for a subset $S \subseteq G$ the
\emph{membership problem for $S$ within $G$}:

\medskip

\noindent
INPUT: A word $w \in (\Sigma^{\pm 1})^*$.

\noindent QUESTION: $h(w) \in S$?

\medskip

\noindent
The \emph{free group} $F(\Sigma)$ generated by $\Sigma$
can be defined as the quotient monoid
$$
F(\Sigma) = (\Sigma^{\pm 1})^*/\{ aa^{-1} = \varepsilon \mid a \in
\Sigma^{\pm 1} \}.
$$
As usual, the \emph{free product} of two groups $G_1$ and $G_2$ is
denoted by $G_1 \ast G_2$. We will always assume that $G_1 \cap G_2 =
\emptyset$.
An \emph{alternating word} in $G_1 \ast G_2$ is a sequence
$g_1 g_2\cdots g_m$ with $m \geq 0$, $g_i \in G_1 \cup G_2$,
and $g_i \in G_1 \Leftrightarrow g_{i+1} \in G_2$. Its length is $m$.
The alternating word $g_1 g_2\cdots g_m$ is \emph{irreducible} if
$g_i \neq 1$ for every $1 \leq i \leq m$.
Every element of $G_1 \ast G_2$ can be written uniquely as an alternating
irreducible word. We will need the following simple fact about
free products:

\begin{lemma}
\label{L free product 1}
Let $g_1 g_2\cdots g_m$ be an alternating word in $G_1\ast G_2$.
If $g_1 g_2\cdots g_m = 1$ in $G_1 \ast G_2$,
then one of the following three cases holds:
\begin{enumerate}[(1)]
\item $m \leq 1$
\item there exists $1 \leq i < m$ such that
$g_1 g_2\cdots g_i = g_{i+1} \cdots g_m = 1$ in  $G_1 \ast G_2$
\item there exist $i \in \{1,2\}$, $k \geq 2$, and
$1 = j_1 < j_2 < \cdots < j_k = m$ such that
$g_{j_1}, g_{j_2},\ldots, g_{j_k} \in G_i$,
$g_{j_1} g_{j_2}\cdots g_{j_k} = 1$ in $G_i$, and
$g_{j_\ell + 1} g_{j_\ell + 2}\cdots g_{j_{\ell+1}-1} = 1$
in  $G_1 \ast G_2$ for all $1 \leq \ell < k$.
\end{enumerate}
\end{lemma}

\begin{figure}[t]
\begin{center}
\setlength{\unitlength}{1mm}
\begin{picture}(86,8)(0,-2)
\drawpolygon[fillgray=1](0,0)(6,0)(6,6)(0,6)
\put(0,0){\makebox(6,6){$g_1$}}
\drawpolygon[fillgray=0.6](6,0)(20,0)(20,6)(6,6)
\put(4.7,0){$\underbrace{\makebox(13.5,6){}}_{\displaystyle 1 \text{ in } G_1\ast G_2}$}
\drawpolygon[fillgray=1](20,0)(26,0)(26,6)(20,6)
\put(20,0){\makebox(6,6){$g_{j_2}$}}
\drawpolygon[fillgray=0.6](26,0)(40,0)(40,6)(26,6)
\put(24.7,0){$\underbrace{\makebox(13.5,6){}}_{\displaystyle 1 \text{ in } G_1\ast G_2}$}
\drawpolygon[fillgray=1](40,0)(46,0)(46,6)(40,6)
\put(40,0){\makebox(6,6){$g_{j_3}$}}
\drawpolygon[fillgray=0.6](46,0)(60,0)(60,6)(46,6)
\put(44.7,0){$\underbrace{\makebox(13.5,6){}}_{\displaystyle 1 \text{ in } G_1\ast G_2}$}
\drawpolygon[fillgray=1](60,0)(66,0)(66,6)(60,6)
\put(60,0){\makebox(6,6){$g_{j_4}$}}
\drawpolygon[fillgray=0.6](66,0)(80,0)(80,6)(66,6)
\put(64.7,0){$\underbrace{\makebox(13.5,6){}}_{\displaystyle 1 \text{ in } G_1\ast G_2}$}
\drawpolygon[fillgray=1](80,0)(86,0)(86,6)(80,6)
\put(80,0){\makebox(6,6){$g_{m}$}}
\end{picture}
\end{center}
\caption{\label{figure case 3}Case (3) in Lemma~\ref{L free product 1}: we have $g_1 g_{j_2}
g_{j_3}g_{j_4}g_m = 1$ in either $G_1$ or $G_2$}
\end{figure}
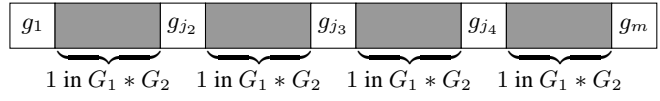

\begin{proof}
Case (3) from the Lemma is visualized in Figure~\ref{figure case 3} for
$k=5$.
Shaded areas represent alternating sequences, which are equal
to $1$ in $G_1 \ast G_2$. The non-shaded blocks are either all
from $G_1$ or from $G_2$, and their product equals $1$ in
$G_1$ or $G_2$, respectively.

We prove the lemma by induction on $m$, the case $m \leq 1$ being trivial.
So assume that $m \geq 2$. Since $g_1 g_2\cdots g_m = 1$ in $G_1
\ast G_2$, there must exist $1 \leq j \leq m$ with $g_j = 1$. If
$j=1$ or $j=m$, then we are in case (2) from the lemma. Hence, we may
assume that $m \geq 3$ and that $2 \leq j \leq m-1$. It follows
$$g_1 \cdots g_{j-2} (g_{j-1}g_{j+1}) g_{j+2} \cdots g_m = 1$$ in
$G_1 \ast G_2$. Since the alternating word $g_1 \cdots g_{j-2}
(g_{j-1}g_{j+1}) g_{j+2} \cdots g_m$ has length $m-2$, we can apply
the induction hypothesis to it. If $m-2 = 1$, i.e., $m=3$, then we
obtain case (3) from the lemma (with $k=2$, $j_1=1$, and $j_2=3$).
If a non-empty and
proper prefix of $g_1 \cdots g_{j-2} (g_{j-1}g_{j+1}) g_{j+2} \cdots
g_m$ equals $1$ in the group $G_1 \ast G_2$, then the same is true
for $g_1 g_2\cdots g_m$. Finally, if case (3) from the lemma applies
to the alternating word $g_1 \cdots g_{j-2} (g_{j-1}g_{j+1}) g_{j+2}
\cdots g_m$, then again the same is true for $g_1 g_2\cdots g_m$.
\qed
\end{proof}
Notice that (3) in Lemma~\ref{L free product 1} can only occur when $m$ is odd.

Assume that $A \leq G$ and $B \leq H$ are groups and 
$\varphi : A \to B$ is an isomorphism. The \emph{amalgamated 
free product} $G\ast_\varphi H$ is the quotient \[(G\ast H)/\{ a = \varphi(a) \mid a \in
A\}.\] Without loss of generality we may assume that $A = G \cap H$ and that 
$\varphi$ is the identity map on $A$; in this situation
we briefly write $G\ast_A H$ for $G\ast_\varphi H$.
Every element of $G\ast_A H$ can be written as a word 
$c_1 \cdots c_n$, where $n \geq 0$, 
$c_1, \ldots, c_n \in G \cup H$, if $n > 1$ then $c_1, \ldots, c_n \in (G \cup H)
\setminus A$, if $n = 1$ then $c_1 \neq 1$, and 
$c_i \in G \setminus A \Leftrightarrow c_{i+1} \in H \setminus A$ for all $1 \leq i
< n$. Such a word is called a \emph{reduced sequence}. 
The normal form theorem for amalgamated 
free products states that every nonempty reduced sequence represents
a nontrivial element of $G\ast_A H$ \cite[Chapter~IV, Theorem~2.6]{LySch77}.

If $G$ is a group and $\varphi:A\to B$ is an isomorphism between
subgroups $A,B$ of $G$, then the \emph{HNN extension} $\bigast_{\varphi} G$,
with base $G$, stable 
letter $t$, and associated subgroups $A,B$ is the quotient group
\[G\ast \langle t\rangle/\{t^{-1}at = \varphi(a)\mid a\in A\}\] where
$t$ is the generator of an infinite cyclic group.  Every element of
$\bigast_{\varphi} G$ can be written as a word  
$g_0t^{\varepsilon_1}g_1\cdots t^{\varepsilon_n}g_n$, where $n \geq 0$, 
$g_0, \ldots, g_n \in G$, and $\varepsilon_1, \ldots, \varepsilon_n \in\{1, -1\}$.
Such a word is referred to as a \emph{reduced sequence} if it contains
no factor of the form $t^{-1}at$ or $tbt^{-1}$ with $a\in A$,
respectively $b\in B$.  Britton's Lemma~\cite[Chapter~IV]{LySch77}
says that if $w=g_0t^{\varepsilon_1}g_1\cdots t^{\varepsilon_n}g_n$ is
a reduced sequence with $n\geq 1$, then $w$ represents a nontrivial
element of $\bigast_{\varphi} G$.

We will also consider fundamental groups of finite graphs of groups,
which is a group theoretic construction generalizing free products,
free products with amalgamation, and HNN-extensions, see e.g.
\cite{Se03}. We omit the quite technical definition. In order to
deal with the rational subset membership problem for graph groups,
free products suffice.

\subsection{Trace monoids and graph groups} \label{SS traces}

In the following we introduce some notions from trace theory, see
\cite{Die90lncs,DieRoz95} for more details.
An \emph{independence alphabet} is just a finite
undirected graph $(\Sigma, I)$ without loops.
Hence, $I \subseteq \Sigma \times \Sigma$ is  an
irreflexive and symmetric relation.
The \emph{trace monoid} $\mathbb{M}(\Sigma,I)$ is defined as the quotient
$$
\mathbb{M}(\Sigma,I) = \Sigma^*/\{ ab = ba \mid (a,b) \in I\}.
$$
Elements of $\mathbb{M}(\Sigma,I)$ are called \emph{traces}.
Note that $\mathbb{M}(\Gamma,J)$ is a submonoid of $\mathbb{M}(\Sigma,I)$
in case $(\Gamma,J)$ is an \emph{induced subgraph} of $(\Sigma,I)$.
The latter means that $\Gamma \subseteq \Sigma$ and 
$J = I \cap (\Gamma \times \Gamma)$.
 
Traces can be represented conveniently by {\em dependence graphs},
which are node-labelled directed acyclic graphs. Let $u = a_1\cdots
a_n$ be a word, where $a_i \in \Sigma$. The vertex set of the
dependence graph of $u$ is $\{1,\ldots,n\}$ and vertex $i$ is
labelled with $a_i \in \Sigma$. There is an edge from vertex $i$ to
$j$ if and only if $i<j$ and $(a_i,a_j) \not\in I$. Then, two words
define the same trace in $\mathbb{M}(\Sigma,I)$ if and only if their
dependence graphs are isomorphic. The set of minimal (resp. maximal)
elements of a trace $t \in \mathbb{M}(\Sigma,I)$ is $\min(t) = \{ a
\in \Sigma \mid \exists u \in  \mathbb{M}(\Sigma,I) : t = au \}$
(resp. $\max(t) = \{ a \in \Sigma \mid \exists u \in
\mathbb{M}(\Sigma,I) : t = ua\}$). A \emph{trace rewriting system}
$R$ over $\mathbb{M}(\Sigma,I)$ is just a finite subset of
$\mathbb{M}(\Sigma,I) \times \mathbb{M}(\Sigma,I)$ \cite{Die90lncs}.
We can define the \emph{one-step rewrite relation} $\to_R \;\subseteq
\mathbb{M}(\Sigma,I) \times \mathbb{M}(\Sigma,I)$ by: $x \to_R y$ if
and only if there are $u,v \in \mathbb{M}(\Sigma,I)$ and $(\ell,r)
\in R$ such that $x = u\ell v$ and $y = u r v$ in $\mathbb{M}(\Sigma,I)$. 
%The notion of a
%\emph{confluent} and \emph{terminating} trace rewriting system is
%defined as for other types of rewriting systems~\cite{Book93}. 
A trace $t$ is \emph{irreducible} with respect to $R$ if there does
not exist a trace $u$ with $t \to_R u$. 
%If $R$ is terminating and
%confluent, then for every trace $t$, there exists a unique
%\emph{normal form} $\NF_R(t)$ such that $t \stackrel{*}{\to}_R
%\NF_R(t)$ and $\NF_R(t)$ is irreducible with respect to $R$.

The \emph{graph group} $\mathbb{G}(\Sigma,I)$ is defined as the quotient
$$
\mathbb{G}(\Sigma,I) = F(\Sigma)/\{ab = ba \mid (a,b) \in I\}.
$$
If $(\Sigma,I)$ is the empty graph, i.e., $\Sigma = \emptyset$, then
we set $\mathbb{M}(\Sigma,I) = \mathbb{G}(\Sigma,I) = 1$ (the
trivial group). Note that $(a,b)\in I$ implies $a^{-1}b = ba^{-1}$
in $\mathbb{G}(\Sigma,I)$. Thus, the graph group
$\mathbb{G}(\Sigma,I)$ can be also defined as the quotient
$$
\mathbb{G}(\Sigma,I) = \mathbb{M}(\Sigma^{\pm 1},I)/\{ aa^{-1} = \varepsilon \mid a \in
\Sigma^{\pm 1} \}.
$$
Here, we implicitly extend $I \subseteq \Sigma\times\Sigma$ to $I \subseteq
\Sigma^{\pm 1} \times \Sigma^{\pm 1}$ by setting
$(a^\alpha, b^\beta) \in I$ if and only if $(a,b) \in I$ for $a,b \in \Sigma$
and $\alpha,\beta\in\{1,-1\}$.
Note that $\mathbb{M}(\Sigma,I)$ is a rational subset of
$\mathbb{G}(\Sigma,I)$.

Define a trace rewriting system $R$ over $\mathbb{M}(\Sigma^{\pm 1},I)$
as follows:
\begin{equation}
\label{system R}
R = \{ (aa^{-1}, \varepsilon) \mid a \in \Sigma^{\pm 1} \}.
\end{equation}
%One can show that $R$ is terminating and confluent and
%that for all $u,v \in \mathbb{M}(\Sigma^{\pm 1},I)$:
One can show that for every trace $t\in \mathbb{M}(\Sigma^{\pm 1},I)$, there exists a unique
\emph{normal form} $\NF_R(t)$ such that $t \stackrel{*}{\to}_R
\NF_R(t)$ and $\NF_R(t)$ is irreducible with respect to $R$.
Moreover, for all $u,v \in \mathbb{M}(\Sigma^{\pm 1},I)$,
$u = v$ in $\mathbb{G}(\Sigma,I)$ if and only if
$\NF_R(u) = \NF_R(v)$ (in $\mathbb{M}(\Sigma^{\pm 1},I)$) \cite{Die90tcs}.
This leads to a linear time solution for the word problem
of $\mathbb{G}(\Sigma,I)$ \cite{Die90tcs,Wra88}.

If the graph $(\Sigma,I)$ is the disjoint union of two graphs
$(\Sigma_1,I_1)$ and $(\Sigma_2,I_2)$, then $\mathbb{G}(\Sigma,I) =
\mathbb{G}(\Sigma_1,I_1) \ast \mathbb{G}(\Sigma_2,I_2)$.  If
$(\Sigma,I)$ is obtained from $(\Sigma_1,I_1)$ and $(\Sigma_2,I_2)$
by connecting each element of $\Sigma_1$ to each element of
$\Sigma_2$, then $\mathbb{G}(\Sigma,I) = \mathbb{G}(\Sigma_1,I_1)
\times \mathbb{G}(\Sigma_2,I_2)$. Graph groups were studied e.g. in
\cite{Dro85}; they are also known as \emph{free partially
commutative groups} \cite{Die90tcs,Wra88}, \emph{right-angled Artin
groups} \cite{BrMe01,CrWi04}, and \emph{semifree groups}
\cite{Bau81}.

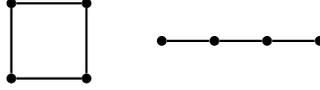
\begin{figure}[t]
\begin{center}
\setlength{\unitlength}{1mm}
\begin{picture}(41,10) 
\gasset{Nframe=n,AHnb=0} 
\gasset{Nw=1.3,Nh=1.3,Nfill=y,linewidth=0.3}
\node(a1)(0,0){}
\node(a2)(10,0){}
\node(a3)(10,10){}
\node(a4)(0,10){}
\drawedge(a1,a2){}
\drawedge(a3,a2){}
\drawedge(a3,a4){}
\drawedge(a1,a4){}

\node(b1)(20,5){}
\node(b2)(27,5){}
\node(b3)(34,5){}
\node(b4)(41,5){}
\drawedge(b1,b2){}
\drawedge(b2,b3){}
\drawedge(b3,b4){}
\end{picture}
\end{center}
\caption{The graphs $\mathsf{C4}$ and $\mathsf{P4}$}
\label{F C4 and P4}
\end{figure}

A \emph{transitive forest} is an independence alphabet $(\Sigma,I)$ such
that there exists a forest $F$ of rooted trees (i.e., a disjoint
union of rooted trees) with node set $\Sigma$ and such that for all
$a, b \in \Sigma$ with $a \neq b$: $(a,b) \in I$ if and only if $a$
and $b$ are comparable in $F$ (i.e., either $a$ is a proper
descendant of $b$ or $b$ is a proper descendant of $a$). It can be
shown that $(\Sigma,I)$ is a transitive forest if and only if
$(\Sigma,I)$ does not contain an induced subgraph, which is a cycle
on 4 nodes (also called $\mathsf{C4}$, see Figure~\ref{F C4 and P4} on the left) or a simple path on 4 nodes
(also called $\mathsf{P4}$, see Figure~\ref{F C4 and P4} on the right) \cite{Wolk65}. The next lemma follows
easily by induction.  We sketch the proof.

\begin{lemma}
\label{L trans. forest}
The class $C$ of all groups, which are of the form
$\mathbb{G}(\Sigma,I)$ for a transitive forest $(\Sigma,I)$, is
the smallest class such that:
\begin{enumerate}[(1)]
\item $1 \in C$
\item if $G_1, G_2 \in C$, then also $G_1 \ast G_2 \in C$
\item if $G \in C$ then $G \times \mathbb{Z} \in C$
\end{enumerate}
\end{lemma}

\begin{proof}
First we verify that graphs groups associated to transitive forests
satisfy (1)-(3).  Case (1) results from the empty graph.
It is immediate that transitive forests are closed under disjoint
union, which implies (2).  If $F$ is a forest of rooted
trees, then one can obtain a rooted tree by adding a new root whose
children are the roots of the trees from $F$.  On the group level
this corresponds to (3).

For the converse, we proceed by induction on the number of vertices.
If the forest $(\Sigma,I)$ consists of more than one rooted tree, then
$\mathbb{G}(\Sigma,I)$ is the free product of the graph groups
associated to the various rooted trees in $(\Sigma,I)$, all of which have a
smaller number of vertices. If there is a single tree, then in
$(\Sigma,I)$ the root is connected to every other vertex. Thus
$\mathbb{G}(\Sigma,I) = G\times \mathbb{Z}$ where $G$ is the graph
group corresponding to the transitive forest obtained by removing
the vertex corresponding to the root and making its children the
roots of the trees in the forest so obtained. \qed\end{proof}

Of course, a similar statement is true for trace monoids
of the form $\mathbb{M}(\Sigma,I)$ with $(\Sigma,I)$
a transitive forest; one just has to replace in
(3) the group $\mathbb{Z}$ by the monoid $\mathbb{N}$.

\section{The rational subset membership problem}
\label{S RSMP}

Let $\mathcal C$ be the smallest class of groups such that:
\begin{itemize}
\item the trivial group $1$ belongs to $\mathcal C$
\item if $G \in \mathcal{C}$ and $H \leq G$ is finitely generated,
  then also $H \in \mathcal{C}$
\item if $G \in \mathcal{C}$ and $G \leq H$ such that $G$ has finite index in $H$
(i.e., $H$ is a finite extension of $G$),
then also $H \in \mathcal{C}$
\item if $G \in \mathcal{C}$, then also $G \times \mathbb{Z} \in \mathcal{C}$
\item if $\mathbb{A}$ is a finite graph of groups~\cite{Se03} whose edge
  groups are finite and whose vertex groups belong to $\mathcal{C}$, then the
  fundamental group of $\mathbb{A}$ belongs to $\mathcal{C}$ (in particular,
  the class $\mathcal C$ is closed under free products).
\end{itemize}
This last property is equivalent to saying that $\mathcal{C}$ is
closed under taking amalgamated products over finite groups and
HNN-extensions with finite associated subgroups~\cite{Se03}. The
main result in this section is:

\begin{theorem}
\label{T dec}
For every group $G \in \mathcal{C}$, the
rational subset membership problem is decidable.
\end{theorem}

It is well known that decidability of the rational subset membership
problem is preserved under taking finitely generated subgroups and finite extensions
\cite{Gru92}. Moreover, the decidability of the rational subset
membership problem is preserved by graph of group constructions with
finite edge groups~\cite{KaSiSt06}.  Hence, in order to prove
Theorem~\ref{T dec}, it would suffice to show that the decidability
of the rational subset membership problem is preserved under direct
products by $\mathbb{Z}$. But currently we can neither prove  nor
disprove this. This forces us to adopt an alternate strategy: we will
introduce an abstract property of groups that implies the decidability of the
rational subset membership problem, and which has the desired
closure properties.

Let $\mathcal L$ be a class of formal languages closed under inverse
homomorphism.  A finitely generated group $G$ is said to be an
\emph{$\mathcal L$-group} if $\text{WP}_{\Sigma}(G)$ belongs to
$\mathcal L$ for some finite generating set $\Sigma$.  This notion is
independent of the choice of generating
set~\cite{Gilman96,HoReRoTh05,KaSiSt06}.

A language $L_0\subseteq \Sigma^*$ belongs to the class \emph{RID}
(rational intersection decidable) if there is an algorithm that,
given a finite automaton over $\Sigma$ recognizing a rational
language $L$, can determine whether $L_0\cap L\neq \emptyset$.  It
was shown in~\cite{KaSiSt06} that the class RID is closed under
inverse homomorphism and that a group $G$ has a decidable rational
subset membership problem if and only if it is an RID-group.  This follows from
the fact that if $L$ is a rational subset of a group $G$, then $g\in
L$ if and only if $1\in g^{-1}L$ and that $g^{-1}L$ is again a
rational subset.

Let $K \subseteq \Theta^*$ be a language over an alphabet $\Theta$.
Then $K$ belongs to the class SLI (semilinear intersection)
if, for every finite alphabet $\Gamma$ (disjoint from $\Theta$) and
every rational language $L \subseteq (\Theta\cup\Gamma)^*$, the set
\begin{equation}\label{setwewantsemilinear}
\Psi_{\Gamma}(\{  w \in L \mid \pi_{\Theta}(w) \in K \}) =
\Psi_\Gamma(L \cap \pi^{-1}_{\Theta}(K))
\end{equation}
is semilinear, and the tuples in a semilinear representation of this
set can be effectively computed from $\Gamma$ and a finite automaton
for $L$. This latter effectiveness statement will be always
satisfied throughout the paper, and we shall not explicitly check
it.  In words, the set \eqref{setwewantsemilinear} is obtained by
first taking those words from $L$ that project into $K$ when
$\Gamma$-letters are erased, and then erasing the $\Theta$-letters,
followed by taking the Parikh image.

In a moment, we shall see that the class SLI is closed under inverse
homomorphism, hence the class of SLI-groups is well defined. In
fact, we show more generally that the class SLI is closed under
inverse images by sequential functions \cite{Ber79}.  This will
imply, moreover, that the class of SLI-groups is closed under taking
finite extensions~\cite{Gilman96,HoReRoTh05,KaSiSt06}.

A {\em sequential transducer} $A$ with input alphabet $\Sigma$ and 
output alphabet $\Omega$ can be defined as a finite state 
automaton with transitions labeled by elements from the set $\Sigma \times \Omega^*$
such that the following restriction is satisfied:
If there are states $q, q_1, q_2$ and a transition from 
$q$ to $q_i$ (for $i \in \{1,2\}$) with label $(a,w_i) \in \Sigma \times \Omega^*$ then 
$q_1 = q_2$ and $w_1 = w_2$. This is not the standard definition
of a sequential transducer (see e.g. \cite{Ber79}), but it
is easily seen to be equivalent. 
The language defined by $A$ is a relation $R \subseteq \Sigma^* \times \Omega^*$,
and it is easy to see that $R$ is the graph of a 
partial function $f:\Sigma^*\to \Omega^*$. A {\em sequential function}
is a partial function, which is computed by a sequential transducer.

\begin{lemma}
\label{L transduction}
Let $K \subseteq \Theta^*$ belong to SLI and let $f:\Sigma^*\to
\Theta^*$ be a sequential function.  Then
$f^{-1}(K)$ belongs to SLI.  In particular, the class of SLI-groups is
well defined and is closed under taking finite extensions.
\end{lemma}

\begin{proof}
Let $\Gamma$ be an alphabet disjoint from $\Sigma$ and
let $L$ be a rational subset of $(\Gamma\cup \Sigma)^*$. Let $A$ be
a sequential transducer computing the sequential function $f : \Sigma^*\to \Omega^*$.
Define a transducer $A'$ by adding to each state of $A$ a loop with label $(a,a)$
for each $a\in \Gamma$. Clearly, $A'$ is a sequential transducer, which
computes a sequential function  
$F:(\Gamma\cup \Sigma)^*\to (\Gamma\cup\Theta)^*$.

The following two observations are immediate from the fact that the only
transitions of $A'$ involving letters from $\Gamma$ are loops
with labels of the form $(a,a)$:
\begin{enumerate}[(a)]
\item $\Psi_{\Gamma} F$ coincides with $\Psi_{\Gamma}$  on the domain of $F$\label{enum a} (we read the composition of
  functions from right to left, i.e., in $\Psi_{\Gamma} F$ we first apply $F$, followed by $\Psi_\Gamma$)
\item $\pi_{\Theta} F = f\pi_{\Sigma}$. \label{enum b}
\end{enumerate}
We now claim that the following equality holds:
\begin{equation}\label{beneq1}
F(L\cap \pi_{\Sigma}^{-1}(f^{-1}(K))) = F(L)\cap \pi^{-1}_{\Theta}(K).
\end{equation}
First note that $L\cap \pi_{\Sigma}^{-1}(f^{-1}(K)) = L\cap
F^{-1}(\pi^{-1}_{\Theta}(K))$ by \eqref{enum b}.  So if $w$ belongs to
the left hand side of \eqref{beneq1}, then $w=F(u)$ with $u\in L\cap
F^{-1}(\pi^{-1}_{\Theta}(K))$.  Thus $w\in F(L)\cap
\pi^{-1}_{\Theta}(K)$.  Conversely, if $u\in F(L)\cap
\pi^{-1}_{\Theta}(K)$, then there exists $w\in L$ such that
$F(w)=u$. But then $w\in L\cap F^{-1}(\pi^{-1}_{\Theta}(K)) =
L\cap \pi_{\Sigma}^{-1}(f^{-1}(K))$ and so
$u$ belongs to the left hand side of \eqref{beneq1}.

Now, since $L\cap \pi_{\Sigma}^{-1}(f^{-1}(K)) = L\cap
F^{-1}(\pi^{-1}_{\Theta}(K))$ is contained in the domain of $F$, we
may conclude from \eqref{enum a} and \eqref{beneq1} that
\begin{equation}\label{beneq2}
\Psi_{\Gamma}(L\cap
\pi_{\Sigma}^{-1}(f^{-1}(K)))=
\Psi_{\Gamma}F(L\cap \pi_{\Sigma}^{-1}(f^{-1}(K))) =
\Psi_{\Gamma}(F(L)\cap \pi^{-1}_{\Theta}(K)).
\end{equation}
But $F(L)$ is rational since the class of rational languages is
closed under images via sequential
functions~\cite{Ber79}. Therefore, since $K$ belongs to SLI, we may
deduce that the Parikh-image $\Psi_{\Gamma}(F(L)\cap \pi^{-1}_{\Theta}(K))$ is
semilinear. This completes the proof of the first statement from the
lemma in light on \eqref{beneq2}.

Since a homomorphism is a sequential function, the
language class SLI is closed under inverse homomorphism. Hence, the
class of SLI-groups is well defined. Finally, let us assume that
$G$ is an SLI-group and that $G$ is a finite index subgroup of $H$.
Let $\Sigma$ (resp. $\Delta$) be a finite generating set for $G$
(resp. $H$). Then in \cite[Lemma 3.3]{KaSiSt06} it is shown that
there exists a sequential function
 $f : \Delta^* \to \Sigma^*$ such that $\text{WP}_{\Delta}(H) =
f^{-1}(\text{WP}_{\Sigma}(G))$. Hence, $H$ is an SLI-group. \qed
\end{proof}

Let us quickly dispense with the decidability of the rational subset
membership problem for SLI-groups.

\begin{lemma}
\label{L 0}
The class of languages SLI is contained in the class of languages RID.
In particular, every SLI-group has a decidable rational subset membership
problem.
\end{lemma}

\begin{proof}
Let $K \subseteq \Theta^*$ belong to SLI.
Let $A$ be a finite automaton over the alphabet $\Theta$. We have
to decide whether $L(A)\cap K\neq \emptyset$. Since $K$ belongs to SLI, the set
$$
\Psi_{\emptyset}(\{  w \in L(A) \mid
\pi_{\Theta}(w) \in
 K \}) =
\Psi_{\emptyset}(L(A) \cap K)
$$
is effectively semilinear and so has a decidable membership problem
(c.f.~\cite{KaSiSt06}). As mentioned earlier,
$\Psi_{\emptyset}(L(A)\cap K)$ consists of the unique function
$\emptyset\to \mathbb{N}$ if $L(A)\cap K$ is non-empty and is empty
otherwise.  Thus we can test emptiness for $L(A)\cap K$.
\qed
\end{proof}

Having already taken care of finite extensions by Lemma~\ref{L transduction},
let's turn to finitely
generated subgroups. We show that the language class SLI is closed
under intersection with rational subsets.
This guarantees that the class of SLI-groups is
closed under taking finitely generated
subgroups~\cite{HoReRoTh05}.

\begin{lemma}
\label{L subgroup} Let $K\subseteq \Theta^*$ belong to SLI and
let $R\subseteq \Theta^*$ be rational.  Then $R\cap K$ belongs
to SLI.  In particular, every finitely generated subgroup of an
SLI-group is an SLI-group.
\end{lemma}

\begin{proof}
Let $L\subseteq (\Gamma\cup \Theta)^*$ be rational, where
$\Gamma$ is a finite alphabet disjoint from $\Theta$.  We have
$$L\cap \pi_{\Theta}^{-1}(R\cap K) = L\cap \pi_{\Theta}^{-1}(R)\cap
\pi_{\Theta}^{-1}(K).$$  
But rational languages are closed under
inverse homomorphism and intersection, so $\Psi_{\Gamma}(L\cap
\pi_{\Theta}^{-1}(R)\cap \pi_{\Theta}^{-1}(K))$ is semilinear as $K$
belongs to SLI.  This establishes the lemma.
\qed
\end{proof}

Next, we show that the class of SLI-groups is closed under
direct products with $\mathbb{Z}$:

\begin{lemma}
\label{L direct product} If $G$ is an SLI-group, then
$G \times \mathbb{Z}$ is also an SLI-group.
\end{lemma}

\begin{proof}
Let $\Sigma$ be a finite generating set for $G$.
Choose a generator
$a\not\in\Sigma$ of $\mathbb{Z}$. Then $G\times\mathbb{Z}$ is generated
by $\Sigma \cup \{a\}$. Let $\Gamma$ be a finite alphabet ($\Gamma \cap
(\Sigma^{\pm 1} \cup \{a,a^{-1}\}) = \emptyset$) and
let $L$ be a rational
subset of $(\Sigma^{\pm 1} \cup \{a,a^{-1}\} \cup \Gamma)^*$.
We have
\begin{align*}
& \Psi_{\Gamma}\Bigl(\{ w \in L \mid
\pi_{\Sigma^{\pm 1}\cup \{a,a^{-1}\}}(w) \in
 \text{WP}_{\Sigma\cup\{a\}}(G \times \mathbb{Z}) \}\Bigr)
= \\
& \qquad\qquad \overline{\pi}_{\Gamma} \Bigl( \; \Psi_{\Gamma\cup
\{a,a^{-1}\}}(\{  w \in L \mid \pi_{\Sigma^{\pm 1}}(w) \in
 \text{WP}_{\Sigma}(G) \})\; \cap \\
& \qquad\qquad \hspace{.67cm}
\{ f \in \mathbb{N}^{\Gamma\cup \{a,a^{-1}\}} \mid f(a)=f(a^{-1}) \} \Bigl).
\end{align*}
This set is semilinear, since $\{ f \in \mathbb{N}^{\Gamma\cup
\{a,a^{-1}\}} \mid f(a)=f(a^{-1}) \}$ is semilinear and semilinear
sets are closed under intersection and projection \cite{GS64}. \qed
\end{proof}

By Lemmas~\ref{L transduction}--\ref{L direct product},
Theorem~\ref{T dec} would be established, if we could prove
the closure of $\mathcal C$ under graph of groups constructions with finite edge
groups. Unfortunately we are only able to prove this closure
under the restriction that every vertex group of the graph of
groups is residually finite (which is the case for groups in $\mathcal C$).
In general we can just prove closure under free product.  This, in
fact, constitutes the most difficult part of the proof of
Theorem~\ref{T dec}.

\begin{lemma}
\label{L free product} If $G_1$ and $G_2$ are
SLI-groups, then  $G_1 \ast G_2$ is also an SLI-group.
\end{lemma}

\begin{proof}
Assume that $\Sigma_i$ is a finite
generating set for $G_i$.
Thus, $\Sigma = \Sigma_1 \cup \Sigma_2$ is
a generating set for the free product $G_1 \ast G_2$.
Let $\Gamma$ be a finite alphabet ($\Gamma \cap
\Sigma^{\pm 1} = \emptyset$) and let
$\Theta = \Sigma^{\pm 1}\cup\Gamma$.
Let $L \subseteq \Theta^*$
be rational and let $A = (Q, \Theta, \delta, q_0, F)$ be a finite automaton
with $L=L(A)$, where $Q$ is the set of states,
$\delta \subseteq Q \times \Theta \times Q$
is the transition relation, $q_0 \in Q$ is the initial state, and
$F \subseteq Q$ is the set of final states.
For $p,q \in Q$ and $w \in \Theta^*$ we write
$p \xrightarrow{w}_A q$ if there exists a path
in $A$ from $p$ to $q$, labelled by the word $w$.

For every pair of states $(p,q) \in Q \times Q$
let us define the language
$$L[p,q]  \subseteq
( \Sigma_1^{\pm 1} \cup \Gamma \cup (Q \times Q))^* \cup
( \Sigma_2^{\pm 1} \cup \Gamma \cup (Q \times Q))^*
\subseteq  (\Theta \cup (Q \times Q))^*$$  as follows:
\begin{align*}
& L[p,q] = \bigcup_{i\in\{1,2\}}\{ w_0 (p_1,q_1) w_1 (p_2,q_2) \cdots w_{k-1} (p_k,q_k) w_k \mid  \\
& \qquad\qquad\qquad k \geq 1 \;\wedge\; (p_1,q_1), \ldots, (p_k,q_k) \in Q\times Q \;\wedge \\
& \qquad\qquad\qquad  w_0, \ldots, w_k \in (
\Sigma_i^{\pm 1} \cup \Gamma )^* \;\wedge\;
\pi_{\Sigma_i^{\pm 1}}(w_0 \cdots w_k)  \in
\text{WP}_{\Sigma_i}(G_i)  \;\wedge \\
& \qquad\qquad\qquad p \xrightarrow{w_0}_A p_1  \;\wedge\; q_1 \xrightarrow{w_1}_A p_2
\;\wedge\cdots\wedge\; q_{k-1} \xrightarrow{w_{k-1}}_A p_k \;\wedge\;
q_k \xrightarrow{w_k}_A q \}
\end{align*}
Since the language
\begin{align*}
& \{ w_0 (p_1,q_1) w_1 (p_2,q_2) \cdots w_{k-1} (p_k,q_k) w_k \mid  \\
& \qquad k \geq 1 \;\wedge\; (p_1,q_1), \ldots, (p_k,q_k) \in Q\times Q \;\wedge \; w_0, \ldots, w_k \in (
\Sigma_i^{\pm 1} \cup \Gamma )^* \;\wedge \\
& \qquad p \xrightarrow{w_0}_A p_1  \;\wedge\; q_1 \xrightarrow{w_1}_A p_2
\;\wedge\cdots\wedge\; q_{k-1} \xrightarrow{w_{k-1}}_A p_k \;\wedge\;
q_k \xrightarrow{w_k}_A q \}
\end{align*}
is a rational language over the alphabet $\Sigma_i^{\pm 1} \cup \Gamma \cup (Q \times Q)$
for $i \in \{1, 2 \}$ and $G_i$ is an SLI-group,
it follows that the Parikh image $\Psi_{\Gamma \cup (Q\times Q)}(L[p,q])
\subseteq \mathbb{N}^{\Gamma \cup (Q \times Q)}$
is semilinear. Let $K[p,q] \subseteq (\Gamma \cup (Q\times Q))^*$ be
some rational language such that
\begin{equation}\label{G K[p,q]}
\Psi(K[p,q]) = \Psi_{\Gamma \cup (Q\times Q)}(L[p,q]).
\end{equation}
Next, we define a context-free grammar $\mathsf{G} = (N,\Gamma,S,P)$
as follows:
\begin{itemize}
\item the set of nonterminals is $N = \{S\} \cup (Q \times Q)$, where $S$ is a new
symbol not contained in $Q \times Q$.
\item $S$ is the start nonterminal.
\item $P$ consists of the following productions:
\begin{alignat*}{2}
S     & \to (q_0, q_f)        && \text{ for all $q_f \in F$}\\
(p,q) & \to K[p,q]            && \text{ for all $p,q \in Q$} \\
(q,q) & \to \varepsilon       && \text{ for all $q \in Q$}
\end{alignat*}
\end{itemize}
By Parikh's theorem, the Parikh image $\Psi(L(\mathsf{G})) \subseteq
\mathbb{N}^{\Gamma}$ is semilinear.
Thus, the following claim proves the lemma:

\medskip

\noindent
{\em Claim 1.}
$\;\Psi(L(\mathsf{G})) = \Psi_{\Gamma}( \{ w \in L(A) \mid
                   \pi_{\Sigma^{\pm 1}}(w) \in \text{WP}_{\Sigma}(
                   G_1\ast G_2) \})$

\medskip

\noindent
{\em Proof of Claim 1.} We prove the following more general identity for all
$(p,q) \in Q \times Q$:
$$
\Psi(L(\mathsf{G},(p,q))) = \Psi_{\Gamma}(\{ w \in \Theta^* \mid p \xrightarrow{w}_A q
\;\wedge\; \pi_{\Sigma^{\pm 1}}(w) \in \text{WP}_{\Sigma}(
                   G_1\ast G_2) \})
$$
For the inclusion 
\begin{equation} \label{Inclusion 1}
\Psi(L(\mathsf{G},(p,q))) \subseteq \Psi_{\Gamma}(\{ w \in \Theta^* \mid p \xrightarrow{w}_A q
\;\wedge\; \pi_{\Sigma^{\pm 1}}(w)
\in \text{WP}_{\Sigma}(G_1\ast G_2)\})
\end{equation}
assume that $(p,q) \stackrel{*}{\Rightarrow}_{\mathsf{G}} u \in \Gamma^*$.
We show by induction on the length of the $\mathsf{G}$-derivation
$(p,q) \stackrel{*}{\Rightarrow}_{\mathsf{G}} u$ that there exists
a word $w \in \Theta^*$ such that
$p \xrightarrow{w}_A q$, $\pi_{\Sigma^{\pm 1}}(w) \in \text{WP}_{\Sigma}(G_1\ast G_2)$,
and $\Psi(u) = \Psi_{\Gamma}(w)$.

\smallskip

\noindent
{\em Case 1.} $p=q$ and $u = \varepsilon$:
We can choose $w = \varepsilon$.

\smallskip

\noindent
{\em Case 2.} $(p,q) \Rightarrow_{\mathsf{G}} u' \stackrel{*}{\Rightarrow}_{\mathsf{G}} u$ for some
$u' \in K[p,q]$. By (\ref{G K[p,q]}), there exists
a word $v \in L[p,q]$ such that
$\Psi(u') = \Psi_{\Gamma \cup (Q\times Q)}(v)$.
Since $v \in L[p,q]$, there exist $k \geq 1$, $(p_1,q_1), \ldots, (p_k,q_k) \in Q\times Q$,
$i \in \{1,2\}$, and $v_0, \ldots, v_k \in (\Sigma_i^{\pm 1} \cup \Gamma )^*$
such that
\begin{itemize}
\item $p \xrightarrow{v_0}_A p_1, \; q_1 \xrightarrow{v_1}_A p_2, \ldots,
q_{k-1} \xrightarrow{v_{k-1}}_A p_k, \; q_k \xrightarrow{v_k}_A q,$
\item $v = v_0 (p_1,q_1) v_1 (p_2,q_2) \cdots v_{k-1} (p_k,q_k) v_k$, and
\item $\pi_{\Sigma_i^{\pm 1}}(v_0 \cdots v_k)  \in \text{WP}_{\Sigma_i}(G_i)$.
\end{itemize}
Since $u' \stackrel{*}{\Rightarrow}_{\mathsf{G}} u \in \Gamma^*$
and $\Psi(u') = \Psi_{\Gamma \cup (Q\times Q)}(v)$,
there must exist $u_1, \ldots, u_k \in \Gamma^*$
such that
$$
(p_i,q_i) \stackrel{*}{\Rightarrow}_{\mathsf{G}} u_i \quad \text{and} \quad
\Psi(u) = \Psi_\Gamma(v_0) + \cdots + \Psi_\Gamma(v_k) + \Psi(u_1) + \cdots + \Psi(u_k)
$$
for all $1 \leq i \leq k$.
By induction, we obtain words $w_1, \ldots, w_k \in\Theta^*$ such that
for all $1 \leq i \leq k$:
\begin{itemize}
\item $p_i \xrightarrow{w_i}_A q_i$
\item $\pi_{\Sigma^{\pm 1}}(w_i) \in \text{WP}_{\Sigma}(G_1\ast G_2)$, and
\item $\Psi(u_i) = \Psi_{\Gamma}(w_i)$.
\end{itemize}
Let us set
$w = v_0 w_1 v_1 \cdots w_k v_k \in \Theta^*$.
We have:
\begin{itemize}
\item
$p \xrightarrow{v_0}_A p_1 \xrightarrow{w_1}_A q_1 \xrightarrow{v_1}_A p_2 \cdots
p_k \xrightarrow{w_k}_A q_k \xrightarrow{v_k}_A q$, i.e., $p \xrightarrow{w}_A q$,
\item
$\pi_{\Sigma^{\pm 1}}(w) \in \text{WP}_{\Sigma}(G_1\ast G_2)$, and
\item
$\Psi(u) =  \Psi_\Gamma(v_0) + \cdots + \Psi_\Gamma(v_k) +
\Psi(u_1) + \cdots + \Psi(u_k) = \Psi_\Gamma(v_0) + \cdots + \Psi_\Gamma(v_k) +
\Psi_\Gamma(w_1) + \cdots + \Psi_\Gamma(w_k) = \Psi_\Gamma(w)$.
\end{itemize}
This concludes the proof of inclusion (\ref{Inclusion 1}).
For the other inclusion, assume that
$$
p \xrightarrow{w}_A q \quad \text{and}\quad \pi_{\Sigma^{\pm 1}}(w)
\in \text{WP}_{\Sigma}(G_1\ast G_2)
$$ for a word $w \in \Theta^*$.
By induction over the length of the word $w$ we show that
$\Psi_{\Gamma}(w) \in \Psi(L(\mathsf{G},(p,q)))$.

We will make a case distinction according to the three cases
in Lemma~\ref{L free product 1}.
Note that we either have $w \in \Gamma^*$ or
the word $w \in \Theta^*$
can be (not necessarily uniquely) written as
$w = w_1 \cdots w_n$ with $n \geq 1$ such that
$w_i \in ((\Gamma \cup \Sigma_1^{\pm 1})^* \cup
(\Gamma \cup \Sigma_2^{\pm 1})^*) \setminus \Gamma^*$
and $w_i \in (\Gamma \cup \Sigma_1^{\pm 1})^*
\Leftrightarrow w_{i+1} \in (\Gamma \cup \Sigma_2^{\pm 1})^*$.

\smallskip

\noindent
{\em Case 1.} $w \in (\Gamma \cup \Sigma_1^{\pm 1})^*$ (the case
$w \in (\Gamma \cup \Sigma_2^{\pm 1})^*$ is analogous):
Then $\pi_{\Sigma_1^{\pm 1}}(w) \in \text{WP}_{\Sigma_1}(G_1)$.
Together with $p \xrightarrow{w}_A q$, we obtain
$w(q,q) = w(q,q) \varepsilon \in L[p,q]$. Since $(p,q) \to K[p,q]$ and
$(q,q) \to \varepsilon$ are productions of $\mathsf{G}$, there exists a word
$u \in \Gamma^*$ such that $(p,q) \stackrel{*}{\Rightarrow}_{\mathsf{G}} u$
and $\Psi(u) = \Psi_\Gamma(w)$, i.e.,
$\Psi_\Gamma(w) \in \Psi(L(\mathsf{G},(p,q)))$.

\smallskip

\noindent
{\em Case 2.} $w = w_1 w_2$ with $w_1 \neq \varepsilon \neq w_2$ and
$\pi_{\Sigma^{\pm 1}}(w_1), \pi_{\Sigma^{\pm 1}}(w_2)
\in \text{WP}_{\Sigma}(G_1\ast G_2)$.
Then there exists a state $r \in Q$ such that
$$
p \xrightarrow{w_1}_A r \xrightarrow{w_2}_A q.
$$
By induction, we obtain
\begin{align*}
\Psi_{\Gamma}(w_1) & \in \Psi(L(\mathsf{G},(p,r))) \text{ and}\\
\Psi_{\Gamma}(w_2) & \in \Psi(L(\mathsf{G},(r,q))).
\end{align*}
Hence, we get
\begin{eqnarray*}
\Psi_{\Gamma}(w) &=& \Psi_{\Gamma}(w_1) + \Psi_{\Gamma}(w_2)  \\
& \in & \Psi(L(\mathsf{G},(p,r))) + \Psi(L(\mathsf{G},(r,q))) \\
& \subseteq & \Psi(L(\mathsf{G},(p,q))),
\end{eqnarray*}
where the last inclusion holds, since
$(p,r)(r,q) \in L[p,q]$, and so either
$(p,q) \to (p,r)(r,q)$ or $(p,q) \to (r,q)(p,r)$
is a production of $\mathsf{G}$.

\smallskip

\noindent
{\em Case 3.} $w = v_0 w_1 v_1 \cdots w_k v_k$
such that $k \geq 1$,
\begin{itemize}
\item $\pi_{\Sigma^{\pm 1}}(w_i) \in \text{WP}_{\Sigma}(G_1\ast G_2)$
for all $i \in \{1, \ldots, k\}$, and
\item for some $i \in \{1,2\}$:
$v_0, \ldots, v_k \in (\Gamma \cup \Sigma_i^{\pm 1})^* \setminus \Gamma^*$ and
$\pi_{\Sigma_i^{\pm 1}}(v_0 \cdots v_k) \in
  \text{WP}_{\Sigma_i}(G_i)$.
\end{itemize}
There exist states $p_1,q_1, \ldots, p_k,q_k \in Q$ such that
$$
p \xrightarrow{v_0}_A p_1 \xrightarrow{w_1}_A q_1 \xrightarrow{v_1}_A p_2 \cdots
p_k \xrightarrow{w_k}_A q_k \xrightarrow{v_k}_A q.
$$
By induction, we obtain
\begin{equation} \label{Psi_Gamma in Psi}
\Psi_\Gamma(w_i) \in \Psi(L(\mathsf{G},(p_i,q_i)))
\end{equation}
for all $1 \leq i \leq k$.
Moreover,
from the definition of the language $L[p,q]$
we obtain
$$
v = v_0 (p_1,q_1) v_1 (p_2,q_2) \cdots v_{k-1} (p_k,q_k) v_k \in L[p,q].
$$
Hence, by (\ref{G K[p,q]}) there is a word $u' \in K[p,q]$ such that
$\Psi(u') = \Psi_{\Gamma \cup (Q\times Q)}(v)$
and $(p,q) \to u'$ is a production of $G$.
With (\ref{Psi_Gamma in Psi}) we obtain
$$
(p,q) \Rightarrow_{\mathsf{G}} u' \stackrel{*}{\Rightarrow}_{\mathsf{G}} u
$$
for a word $u \in \Gamma^*$ such that
$$
\Psi(u) = \Psi_\Gamma(v_0)+\cdots+\Psi_\Gamma(v_k)+
\Psi_\Gamma(w_1)+\cdots+\Psi_\Gamma(w_k) = \Psi_\Gamma(w),
 $$
i.e., $\Psi_\Gamma(w) \in \Psi(L(\mathsf{G},(p,q)))$.
This concludes the proof of Claim 1 and hence of the lemma.
\qed
\end{proof}

If we were to weaken the definition of the class $\mathcal C$ by
only requiring closure under free products instead of closure under
finite graphs of groups with finite edge groups, then Lemmas~\ref{L
0}--\ref{L free product} would already imply Theorem~\ref{T dec}. In
fact, this weaker result suffices in order to deal with graph
groups, and readers only interested in graph groups can skip the
following considerations concerning graphs of groups.

To obtain the more general closure result for the class $\mathcal C$
concerning graph of group constructions, we reduce to the case of
free products.  Recall that a group $G$ is \emph{residually finite}
if, for each $g\in G\setminus \{1\}$, there is a finite index normal
subgroup $N$ of $G$ with $g\notin N$. Now we use a standard trick
for graphs of residually finite groups with finite edge groups.

\begin{lemma}
\label{L graph of groups} Let $\mathbb{A}$ be a finite graph of groups
such that the vertex groups are residually finite SLI-groups and the edge
groups are finite.
Then the fundamental group of
$\mathbb{A}$ is an SLI-group.
\end{lemma}

\begin{proof}
Let $G$ be the fundamental group of $\mathbb{A}$.  Then $G$ is
residually finite~\cite[II.2.6 Proposition 12]{Se03}.  Since there are only
finitely many  
edge groups and each edge group is finite, there is a finite index
normal subgroup $N \leq G$ intersecting trivially each edge group,
and hence each conjugate of an edge group.  Thus the finitely
generated subgroup $N \leq G$ acts on the Bass-Serre tree for
$G$~\cite{Se03} with trivial edge stabilizers, forcing $N$ to be a
free product of conjugates of subgroups of the vertex groups of $G$
and a free group~\cite{Se03}. Since $N$ is finitely generated, these
free factors must also be finitely generated. Since every finitely
generated subgroup of an SLI-group is an SLI-group (Lemma~\ref{L
subgroup}) and $\mathbb{Z}$ is an SLI-group (Lemma~\ref{L direct
product}), we may deduce that $N$ is a free product of SLI-groups
and hence is an SLI-group by Lemma~\ref{L free product}. Since $G$
contains $N$ as a finite index subgroup, Lemma~\ref{L transduction}
implies that $G$ is an SLI-group, as required. \qed
\end{proof}

Clearly, the trivial group $1$ is an SLI-group.  Also all the
defining properties of $\mathcal C$ preserve residual finiteness
(the only non-trivial case being the graph of group
constructions~\cite{Se03}). Hence, Lemmas~\ref{L 0}--\ref{L direct product}
and Lemma~\ref{L graph of groups} immediately yield Theorem~\ref{T dec}.

Our main application of Theorem~\ref{T dec} concerns graph groups:

\begin{theorem}
\label{T main}
The rational subset membership problem for a graph
group $\mathbb{G}(\Sigma,I)$ is decidable if and only if
$(\Sigma,I)$ is a transitive forest. Moreover, if $(\Sigma,I)$
is not a transitive forest, then there exists a fixed rational
subset $L$ of $\mathbb{G}(\Sigma,I)$ such that the membership problem
for $L$ within $\mathbb{G}(\Sigma,I)$ is undecidable.
\end{theorem}

\begin{proof}
The decidability part follows immediately from Theorem~\ref{T dec}:
Lemma~\ref{L trans. forest} implies that every graph group
$\mathbb{G}(\Sigma,I)$ with $(\Sigma,I)$ a transitive forest belongs
to the class $\mathcal C$.

Now assume that $(\Sigma, I)$ is not a transitive forest.
By \cite{Wolk65} it suffices to
consider the case that $(\Sigma,I)$ is either a
$\mathsf{C4}$ or a $\mathsf{P4}$. For the case of a $\mathsf{C4}$ we
can use Mihailova's result \cite{Mih66}. Now assume that
$(\Sigma,I)$ is a $\mathsf{P4}$.  
We will reuse a construction by Aalbersberg and Hoogeboom 
\cite{AaHo89}, which is based on \emph{2-counter machines}.
A 2-counter machine is a tuple $C = (Q, \text{Ins}, q_0, q_f)$
where $Q$ is a finite set of states,
$q_0 \in Q$ is the initial state, $q_f \in Q$ is the final state,
and  $\text{Ins} \subseteq Q \times \{i1, i2, d1, d2, z1, z2, p1, p2\} \times Q$ is
the set of instructions. 
The set of configurations of $C$ is $Q \times \mathbb{N} \times \mathbb{N}$.
For two configurations $(q, n_1, n_2), (q', m_1, m_2)$ we write 
$(q, n_1, n_2) \Rightarrow_C (q', m_1,m_2)$ if there exists an instruction $(q, \alpha k, q') \in
\text{Ins}$, so $\alpha \in \{ i, d, z, p \}, k \in \{1,2\}$,
such that $m_{3-k} = n_{3-k}$ and 
one of the following three cases holds:
\begin{itemize}
\item $\alpha = i$ and $m_k = n_k+1$
\item $\alpha = d$ and $m_k = n_k-1$
\item $\alpha = z$ and $m_k = n_k = 0$
\item $\alpha = p$ and $m_k = n_k > 0$
\end{itemize} 
Since Turing machines can be simulated by 2-counter machines \cite{Min61},
it is undecidable whether for a given 2-counter machine $C = (Q, \text{Ins}, q_0, q_f)$ there 
exist $m,n \in \mathbb{N}$ with $(q_0,0,0) \Rightarrow_C^* (q_f,m,n)$.
In \cite{AaHo89}, this problem is reduced to the question, whether 
$L \cap K = \emptyset$ for given  rational trace languages $L, K
\subseteq \mathbb{M}(\Sigma, I)$, where $\Sigma = \{a,b,c,d\}$ and
$I = \{ (a,b), (b,c), (c,d) \}$.
In fact, the language $K$ is
fixed, more precisely 
\begin{eqnarray*}
K & = & ba( d (cb)^+ a)^* d c^* \\
  & = & \{ [a b^{j_0} c^{j_1}  d a b^{j_1} c^{j_2} d \cdots a b^{j_{\ell-1}}
  c^{j_\ell} d]_I \mid \ell \geq 1, j_0 = 1, j_1,\ldots, j_\ell \geq 1 \} .
\end{eqnarray*}
The problem is
that in the construction of \cite{AaHo89} the language $L$ is not fixed
since it depends on the 2-counter machine $C$. 
Aalbersberg and Hoogeboom encode 
the pair of counter values $(m,n) \in \mathbb{N} \times
\mathbb{N}$ by the single number $2^{m} 3 ^{n}$. The language
$L$ is constructed in such a way that
$K \cap L$ contains exactly those traces of the form 
$[a b^{j_0} c^{j_1}  d a b^{j_1} c^{j_2} d \cdots a b^{j_{\ell-1}} c^{j_\ell}
d]_I$, such that $\ell \geq 1$, $j_0 = 1$, and 
there exist states $q_1, \ldots, q_\ell$ and $m_1, n_1, \ldots, m_\ell, n_\ell \in \mathbb{N}$
with $q_\ell = q_f$, $2^{m_i}3^{n_i} = j_i$, and 
$(q_0,0,0) \Rightarrow_C (q_1, m_1, n_1) \Rightarrow_C (q_2, m_2, n_2) \Rightarrow_C
\cdots \Rightarrow_C (q_\ell, m_\ell, n_\ell)$
(note that $j_0 = 1$ encodes the initial counter values $(0,0)$).

In order to construct a fixed rational
subset of $\mathbb{G}(\Sigma,I)$ with an undecidable membership problem,
we start with a fixed (universal) 2-counter machine $C = (Q, \text{Ins}, q_0, q_f)$
such that it is undecidable whether $\exists m', n' \in \mathbb{N}: 
(q_0, m, n) \Rightarrow_C^* (q_f, m', n')$ for given natural numbers $m,n$.
Such a machine $C$ can be obtained by simulating a universal Turing machine.
Let $L \subseteq \mathbb{M}(\Sigma, I)$
be the \emph{fixed} rational trace language constructed by
Aalbersberg and Hoogeboom from $C$, and let us replace the fixed
trace language $K = ba( d (cb)^+ a)^* d c^*$ by the (non-fixed)
language 
\begin{eqnarray*}
K_{m,n} & = &  b^{2^m 3^n} a( d (cb)^+ a)^* d c^* \\
  & = & \{ [a b^{j_0} c^{j_1}  d a b^{j_1} c^{j_2} d \cdots a b^{j_{\ell-1}}
  c^{j_\ell} d]_I \mid \ell \geq 1, j_0 = 2^m3^n, j_1,\ldots, j_\ell \geq 1 \} .
\end{eqnarray*}
Then it is undecidable, whether $K_{m,n} \cap L \neq \emptyset$
for given $m,n \in \mathbb{N}$.
Hence, it is undecidable, whether
$b^{- 2^m 3^n} \in a( d (cb)^+ a)^* d c^* L^{-1}$ in
the graph group $\mathbb{G}(\Sigma, I)$.
Clearly, $a( d (cb)^+ a)^* d c^* L^{-1}$ is a fixed
rational subset of the graph
group $\mathbb{G}(\Sigma, I)$. \qed
\end{proof}

We conclude this section with a further application of Theorem~\ref{T dec}
to \emph{graph products} (which should not be confused
with graphs of groups). A graph product is given by a tuple
$(\Sigma, I, (G_v)_{v\in \Sigma})$, where $(\Sigma,I)$ is an
independence alphabet and $G_v$ is a group, which is associated
with the node $v\in \Sigma$. The group $\mathbb{G}(\Sigma, I, (G_v)_{v\in \Sigma})$
defined by this tuple is the quotient
$$
\mathbb{G}(\Sigma, I, (G_v)_{v\in \Sigma}) = \bigast_{v \in \Sigma} G_v / \{ xy = yx \mid
x \in G_u, y \in G_v, (u,v) \in I \},
$$
i.e., we take the free product
$\bigast_{v \in \Sigma} G_v$ of the groups $G_v$ ($v \in \Sigma$), but
let elements from adjacent groups commute.
Note that $\mathbb{G}(\Sigma, I, (G_v)_{v\in \Sigma})$ is the graph
group $\mathbb{G}(\Sigma, I)$ in the case every $G_v$ is isomorphic to
$\mathbb{Z}$. Graph products were first studied by Green \cite{Gre90}.

\begin{theorem}
If $(\Sigma, I)$ is a transitive forest and every group $G_v$ ($v \in V$)
is finitely generated and virtually Abelian (i.e., has an Abelian subgroup of
finite index), then the rational subset membership problem for
$\mathbb{G}(\Sigma, I, (G_v)_{v\in \Sigma})$ is decidable.
\end{theorem}

\begin{proof}
Assume that the assumptions from the theorem are satisfied.
We show that $\mathbb{G}(\Sigma, I, (G_v)_{v\in \Sigma})$
belongs to the class $\mathcal C$. Since
$(\Sigma,I)$ is a transitive forest, the group
$\mathbb{G}(\Sigma, I, (G_v)_{v\in \Sigma})$ can be built up
from trivial groups using the
following two operations: (i) free products and (ii)
direct products with  finitely generated virtually Abelian groups.
Since the class $\mathcal C$ is closed
under free products, it suffices to prove
that if $G$ belongs to the class $\mathcal C$ and
$H$ is finitely generated virtually Abelian, then
$G \times H$ also belongs to the class $\mathcal C$.
As a finitely generated virtually Abelian group, $H$ is a finite extension of a finite
rank free Abelian group $\mathbb{Z}^n$. By the closure of the class $\mathcal C$
under direct products with $\mathbb{Z}$,
$G \times \mathbb{Z}^n$ belongs to the class $\mathcal C$.
Now, $G \times H$ is a finite extension
of $G \times \mathbb{Z}^n$, proving the theorem,
since  $\mathcal C$ is closed under finite extensions.
\qed
\end{proof}

\section{The submonoid membership problem}\label{S submonoid}

Recall that the submonoid membership problem for a group $G$ asks
whether a given element of $G$ belongs to a given finitely generated
submonoid of $G$. Hence,
there is a trivial reduction from the submonoid membership problem
for $G$ to the rational subset membership problem for $G$. We will
show that for every amalgamated free product $G \ast_A H$ such that:
\begin{enumerate}
\item $A = G \cap H$ is a finite, proper subgroup of $G$ and $H$;
\item there exist $g \in G$, $h \in H$ with $g^{-1}Ag \cap A = 1 = 
h^{-1}Ah \cap A$,
\end{enumerate}
there is in
fact also a reduction in the opposite direction. Similarly, if
$\bigast_{\varphi}G$ is an HNN extension with $\varphi:A\to B$ with 
\begin{enumerate}
\item $A$ is a finite subgroup of $G$;
\item there exists $g\in G$ such that $g^{-1}Ag\cap A=1$ or
  $g^{-1}Ag\cap B =1$ 
\end{enumerate}
then the rational subset problem reduces to the submonoid membership
problem for $\bigast_{\varphi}G$.  We remark that in 2, one could by
symmetry switch the roles of $B$ and $A$.

Using the following lemma, it will suffice to
consider a free product $G \ast F_2$, where $F_2$
is a free group of rank two.

\begin{lemma}\label{L free group of rank 2}
Let $G\ast_A H$ be an amalgamated free product such that
$H \neq A$, $[G:A] \geq 5$, and there exists $h \in H$ 
with $h^{-1} A h \cap A = 1$.
Then $G\ast_A H$ contains as a subgroup the free product $G \ast F_2$
of $G$ with a free group of rank two.
\end{lemma}

\begin{proof}
Since  $[G:A] \geq 5$, we can choose elements
$g_1, g_2, g_3, g_4 \in G \setminus A$ which belong to 
pairwise distinct left $A$-cosets. Moreover, choose an element
$h \in H\setminus A$ with $h^{-1} A h \cap A = 1$. First
we claim that $x = g_1 h g_2^{-1}$ and $y=g_3 h g_4^{-1}$
freely generate a free subgroup of $G\ast_A H$.  
For this, note that $g_i^{-1} g_j \in G \setminus A$ if $i \neq j$.
Thus, every word over $\{x,x^{-1},y, y^{-1}\}$ which does not contain a 
factor from $\{ x x^{-1}, x^{-1}x, y y^{-1}, y^{-1}y\}$
yields a reduced
sequence for the amalgamated product. The normal form theorem for
amalgamated free products~\cite[Chapter~IV, Theorem~2.6]{LySch77}
then implies that $\{ x, y \}$ is 
the base of a free subgroup of $G \ast_A H$.
Hence, the conjugates
$u = hxh^{-1}$ and $v = hyh^{-1}$ also form a base for a free subgroup
of $G \ast_A H$.
Since $h^{-1} A h \cap A = 1$ (and hence if $a\in A$, then
$h^{-1}ah\in H\setminus A$) a word over $G\setminus\{1\} \cup \{ u,
u^{-1},v,v^{-1}\}$, which does not contain a factor
from $(G \setminus\{1\})(G \setminus\{1\}) \cup \{ u u^{-1}, u^{-1}u,
v v^{-1}, v^{-1}v\}$, yields a reduced sequence for the amalgamated
product. Again,  
the normal form theorem for amalgamated free products implies that 
the subgroup of $G\ast_A H$ generated by
$G \cup \{ u, v\}$ is isomorphic to $G\ast F_2$. 
\qed
\end{proof}

We now prove the analogous result for HNN extensions.

\begin{lemma}\label{L free group of rank 2 HNN version}
Let $\bigast_{\varphi} G$ be an HNN extension with stable letter $t$
and finite associated subgroups $A,B$ (so $\varphi:A\to B$) such that
$[G:B] \geq 3$ and there exists $g \in G$ 
with $g^{-1} Ag  \cap A = 1$ or $g^{-1} Ag\cap B=1$.
Then $\bigast_{\varphi}G$ contains as a subgroup the free product $G \ast F_2$
of $G$ with a free group of rank two.
\end{lemma}

\begin{proof}
By Lemma~\ref{L free group of rank 2}, it suffices to show that
$\bigast_{\varphi} G$ contains a subgroup $G\ast \mathbb Z$. 
We may assume that $A \neq 1 \neq B$, because otherwise
$\bigast_{\varphi}G \simeq G * \mathbb{Z}$. Choose
$g_1,g_2\in G\setminus B$ so that $g_1,g_2$ are in different left cosets of
$B$.  Suppose first there exists $g \in G$ with $g^{-1}Ag\cap A=1$ and
set $x=g_1t^{-1}gtg_2^{-1}$.  Since $g\notin A$ (because otherwise
$A=1$) and
$g_2^{-1}g_1\notin B$, one easily deduces that $x^n$ is a reduced
sequence for the HNN extension for all $n>0$ and hence $x$ is of
infinite order by Britton's lemma.  Set $y = t^{-1}gtxt^{-1}g^{-1}t$.
Then $y$ is of infinite order, being a conjugate of $x$.  We claim
that $G$ and $\langle y\rangle$ generate their free product inside of
$\bigast_{\varphi} G$.  We need to show that a word over $G\setminus \{1\}\cup
\{y,y^{-1}\}$ with no factor from $(G \setminus\{1\})(G \setminus\{1\})
\cup \{yy^{-1}, y^{-1}y\}$ results in a reduced sequence for the HNN
extension.  The key point is that if $h\in G\setminus B$, then
$t^{-1}g^{-1}tht^{-1}gt$ is reduced.  On the other hand, if $b\in B\setminus\{1\}$, then
$t^{-1}g^{-1}tbt^{-1}gt = t^{-1}g^{-1}\varphi^{-1}(b)gt^{-1}$, which
is reduced since $g^{-1}Ag\cap A=1$.

Now assume that there exists $g \in G$ with 
$g^{-1}Ag\cap B=1$. The group $A$ must be a proper subgroup of
$G$, because otherwise we have $1 = g^{-1}Ag \cap B = G \cap B = B$.
So choose $g_0\in G\setminus A$ and set
$x=g_1t^{-1}g_0tg_2^{-1}$.  The same argument as above shows that $x$
has infinite order.  Set $y = t^{-1}gt^{-1}xtg^{-1}t$; again $y$ has
infinite order, being a conjugate of $x$.  Again, we claim that $G$
and $\langle y\rangle$ generate their free product in
$\bigast_{\varphi} G$.  Once 
more, we must prove that a word over $G\setminus \{1\}\cup
\{y,y^{-1}\}$ with no factor from $(G \setminus\{1\})(G \setminus\{1\})
\cup \{yy^{-1}, y^{-1}y\}$ yields a reduced sequence for the HNN
extension.  The key point is that if $h\in G\setminus B$, then
$tg^{-1}tht^{-1}gt^{-1}$ is reduced.  On the other hand, if $b\in B\setminus \{1\}$, then
$tg^{-1}tbt^{-1}gt^{-1} = tg^{-1}\varphi^{-1}(b)gt^{-1}$, which
is reduced since $g^{-1}Ag\cap B=1$.
\qed\end{proof}

The following lemma is crucial for us:

\begin{lemma} \label{L submonoid free prod}
{}\
\begin{enumerate}[(1)]
\item Let $G$ and $H$ be finitely generated groups such that
the finite group $A$ is a proper subgroup of both $G$ and $H$
and there exists $h \in H$ with $h^{-1}Ah \cap A = 1$.
Then the rational subset membership problem for
$G$ can be reduced to the submonoid membership problem for
$G \ast_A H$.  

\item If $\varphi:A\to B$ is an isomorphism between finite subgroups
of a finitely generated group $G$ and there exists $g\in G$ with
$g^{-1}Ag\cap A=1$ or 
$g^{-1}Ag\cap B=1$,  then the rational subset membership problem for $G$ can be
reduced to the submonoid membership problem for $\bigast_{\varphi} G$.
\end{enumerate}
\end{lemma}

\begin{remark}
In our proof of Lemma~\ref{L submonoid free prod} we will implicitly
construct Turing machines that carry out the reductions in (1) and (2).
These machines will depend on the element $g$ (and $h$) in (1), respectively (2).
Here one might argue that these elements are not known. But this
is not a real problem, since $g$ and $h$ are fixed elements which
do not depend on the input for the reduction. So there exists a Turing machine
that can do the reduction, although we don't know which Turing machine
if we don't know the elements $g$ and $h$.
\end{remark}

\noindent
{\it Proof of Lemma~\ref{L submonoid free prod}.}
If $G$ is finite, then the rational subset membership problem for
$G$ is decidable, so we may assume without loss of generality
that $G$ is infinite. Since $A$ is finite, we have $[G:A] \geq 5$ in (1),
respectively $[G:B]\geq 3$ in (2).
Then Lemmas~\ref{L free group of rank 2} 
and~\ref{L free group of rank 2 HNN version} imply that
$G\ast F_2$ is a subgroup of $G \ast_A H$, respectively
$\bigast_{\varphi} G$. Since the submonoid membership
problem for a finitely generated subgroup of a group 
$K$ reduces to the  submonoid membership 
problem for $K$ itself, it suffices to 
prove the following:  the rational subset membership problem for $G$ can be
reduced to the submonoid membership problem for $G \ast F_2$.
%(note that the proofs of Lemmas~\ref{L free group of rank 2} 
%and~\ref{L free group of rank 2 HNN version} are effective).
Let $\Sigma$ be a finite generating set for $G$ and use
$h:(\Sigma^{\pm 1}\cup \Gamma^{\pm})^*\to G\ast F_2$ for the
canonical morphism. Let $A = (Q, \Sigma^{\pm 1}, \delta, q_0, F)$ be
a finite automaton and let $t \in (\Sigma^{\pm 1})^*$. By
introducing $\varepsilon$-transitions, we may assume that the set of
final states $F$ consists of a single state $q_f \neq q_0$.  One can
effectively find a subset $\widetilde Q\subseteq F_2$ in bijection
with $Q$ via $q\mapsto \widetilde q$ such that $\widetilde Q$ freely
generates a free subgroup of $F_2$.

We construct a finite subset $\Delta \subseteq (\Sigma^{\pm 1}\cup
\Gamma^{\pm 1})^*$ and an element $u \in (\Sigma^{\pm 1}\cup
\Gamma^{\pm 1})^*$ such that $h(t) \in h(L(A))$ if and only if
$h(u)\in h(\Delta^*)$. Let
\begin{equation}\label{defineDelta}
\Delta = \{ \widetilde{q}\, c\, \widetilde{p}^{-1} \mid (q, c, p)
\in \delta \} \quad\text{ and }\quad u = \widetilde{q}_0\, t
\,\widetilde{q}_f^{-1} .
\end{equation}
Note that in \eqref{defineDelta}, we have $c \in \Sigma^{\pm 1} \cup
\{ \varepsilon \}$, since we introduced $\varepsilon$-transitions.
Recall $(q,c,p)\in \delta$ means $q\xrightarrow{c}p$ in $A$.
We begin with a critical claim.

\medskip
\noindent {\em Claim 1.}  Suppose that in $G\ast F_2$, we have
\begin{equation}\label{claim1eq}
\widetilde{q}_0\, t\,\widetilde{q}_f^{-1} =
(\widetilde{p}_1 v_1 \widetilde{q}_1^{-1})\cdots
(\widetilde{p}_n v_n \widetilde{q}_n^{-1})
\end{equation}
 where
$p_i\xrightarrow{v_i}q_i$ in $A$, for $i \in \{1,\ldots,n\}$.  
Then $h(t)\in h(L(A))$.
\medskip

\noindent The claim is proved by induction on $n$.  If $n=1$, then
since $q_0\neq q_f$, the  normal form theorem for free products
easily implies $q_0=p_1$, $q_f=q_1$ and $t=v_1$ in $G$.  Thus
$q_0\xrightarrow{v_1}q_f$ in $A$, whence $v_1\in L(A)$, and so
$h(t)\in h(L(A))$.  Next suppose the claim holds for $n-1\geq 1$ and
consider the claim for $n>1$.

First suppose that $q_i = p_{i+1}$ for some $i$.  Then
\[\widetilde{q}_0\, t\,\widetilde{q}_f^{-1} =
(\widetilde{p}_1 v_1 \widetilde{q}_1^{-1})\cdots
(\widetilde{p}_iv_iv_{i+1}\widetilde{q}_{i+1}^{-1})\cdots
(\widetilde{p}_n v_n \widetilde{q}_n^{-1})\] in $G \ast F_2$ and
$p_i\xrightarrow{v_iv_{i+1}}q_{i+1}$ in $A$.  Induction now gives
the desired conclusion.

Next suppose that for some $i$, we have $p_i=q_i$ and $v_i=1$ in $G$.
Then \[\widetilde{q}_0\, t \,\widetilde{q}_f^{-1} =
(\widetilde{p}_1 v_1 \widetilde{q}_1^{-1})\cdots
(\widetilde{p}_{i-1} v_{i-1}
\widetilde{q}_{i-1}^{-1})(\widetilde{p}_{i+1}v_{i+1}\widetilde{q}_{i+1}^{-1})\cdots 
(\widetilde{p}_n v_n \widetilde{q}_n^{-1})\] in $G \ast F_2$ and we
can again apply 
the induction hypothesis.

Finally, suppose  $p_i=q_i$ implies $v_i\neq 1$ in $G$ and suppose
$q_i\neq p_{i+1}$, all $i$.  Then we claim that the right hand side
of \eqref{claim1eq} is already in normal form.   Consider a typical window
$\widetilde{q}_{i-1}^{-1}\widetilde{p}_i v_i
\widetilde{q}_i^{-1}\widetilde{p}_{i+1}$ (where we take
$\widetilde{q}_0=1=\widetilde{p}_{n+1}$). 
Then no two neighbouring elements belong to the same factor of the
free product $G\ast \langle \tilde Q\rangle = G\ast \langle
\widetilde{s}_1\rangle\ast \cdots \ast \langle
\widetilde{s}_m\rangle$, where $Q=\{s_1,\ldots, s_m\}$, since
$q_j\neq p_{j+1}$ for $j=i-1,i$ and $p_i\neq q_i$ when $v_i=1$ in $G$. Since
such windows cover the right hand side of \eqref{claim1eq} we may
conclude that it is in normal form in $G\ast F_2$. Comparison with
the left hand side then shows that $n=1$, contradicting $n>1$.  So
this case does not arise and the proof of the claim is complete.

Now we may prove that $h(t)\in h(L(A))$ if and only if $h(u)\in
h(\Delta^*)$.  Suppose first that $h(t) = h(t')$ with $t'\in L(A)$.
Write $t'=a_1\cdots a_n$ with $a_i\in \Sigma^{\pm 1} \cup
\{\varepsilon\}$ and such 
that $q_0\xrightarrow{a_1}q_1\xrightarrow{a_2}q_2\longrightarrow
\cdots\longrightarrow q_{n-1}\xrightarrow{a_n}q_f$.  Then, as
$h(t)=h(t')$,  clearly we
have
\[ u = \widetilde{q_0} t' \widetilde{q_f}^{-1}=
(\widetilde{q}_0 a_1 \widetilde{q}_1^{-1})(\widetilde{q}_1 a_2 \widetilde{q}_2^{-1})
\cdots (\widetilde{q}_{n-1} a_n \widetilde{q}_f^{-1})\in
\Delta^* \] 
in $G \ast F_2$.
Conversely, suppose $h(u)\in h(\Delta^*)$. Then we
can write \[ u = \widetilde{q}_0\, t\,\widetilde{q}_f^{-1} =
(\widetilde{p}_1 a_1 \widetilde{q}_1^{-1})\cdots
(\widetilde{p}_n a_n \widetilde{q}_n^{-1})\] 
in $G \ast F_2$, where
$p_i\xrightarrow{a_i}q_i$ are certain transitions of $A$.   Claim 1
then implies $h(t)\in h(L(A))$.
\qed

\begin{theorem} \label{T submonoid}
Let $G$ and $H$ be finitely generated groups
such that the finite group $A$ is a proper subgroup of both $G$ and $H$
and there exist $g \in G$, $h \in H$ with $g^{-1}Ag \cap A = 1 = 
h^{-1}Ah \cap A$.
Then, for the amalgamated free product $G \ast_A H$ the
rational subset membership problem and
the submonoid membership problem are recursively
equivalent.
\end{theorem}

\begin{proof}
It suffices to show that the rational subset membership problem
for $G \ast_A H$ can be reduced to the submonoid membership problem
for $G \ast_A H$.
The rational subset membership problem for $G \ast_A H$ can be reduced to
the rational subset membership problems for $G$ and $H$ \cite{KaSiSt06}.
By Lemma~\ref{L submonoid free prod}
both these problems can be reduced to the submonoid membership problem
for $G \ast_A H$.
\qed
\end{proof}

Note that the assumptions in  Theorem~\ref{T submonoid} are satisfied
for every free product $G \ast H$ of nontrivial finitely generated groups
$G$ and $H$.

A similar result holds for HNN extensions:

\begin{theorem} \label{T submonoid HNN}
Let $G$ be a finitely generated group and let $\varphi:A\to B$ be an
isomorphism between finite subgroups of $G$.  Suppose there exists $g
\in G$, with $g^{-1}Ag \cap A = 1$ or $g^{-1}Ag \cap B = 1$.
Then the
rational subset membership problem and
the submonoid membership problem are recursively
equivalent for the HNN extension $\bigast_{\varphi} G$.
\end{theorem}

\begin{proof}
We just need to establish that the rational subset membership problem
for $\bigast_{\varphi} G$ can be reduced to the submonoid membership problem.
The rational subset membership problem for $\bigast_{\varphi} G$ can
be reduced to the rational subset membership problem for $G$ by the
results of~\cite{KaSiSt06}.
By Lemma~\ref{L submonoid free prod} this problem can be reduced to
the submonoid membership problem for $\bigast_{\varphi} G$.  This
completes the proof.
\qed
\end{proof}

Let us say that a group $G$ is \emph{virtually a free product} if it
has a finite index subgroup $H$ that splits nontrivially as a free
product $H=G_1\ast G_2$.

\begin{corollary} 
Let $G$ be a finitely generated group that is virtually a free
product.  Then the rational subset and submonoid membership problems
are recursively equivalent.
\end{corollary}

\begin{proof}
Suppose $G$ has decidable submonoid membership problem.  We
need to show that $G$ has decidable rational subset problem.  Let
$H$ be a finite index subgroup of $G$ that splits nontrivially as a
free product.  Clearly $H$ has decidable submonoid membership
problem and hence has decidable rational subset membership problem
by Theorem~\ref{T submonoid}.  It then follows $G$ has decidable
rational subset membership problem by~\cite{Gru92,KaSiSt06}.
\qed
\end{proof}

In order for a finitely generated group to be virtually a free
product, it must have two or more ends.  On the other hand, a group
with two or more ends that is either virtually torsion-free or
residually finite is easily seen, via Stallings ends
theorem~\cite{St71}, to be virtually a free product, as we now show.
First we recall the notion of  ends of a locally finite graph.

Let $\Gamma$ be a locally finite graph, i.e., every node of $\Gamma$
has only finitely many neighbours. Consider the inverse system
$\Gamma\setminus C$ where $C$ runs over the finite subgraphs of
$\Gamma$.  Then the sets of connected components
$\pi_0(\Gamma\setminus C)$ form an inverse system of sets; the
projective limit $\mathrm{Ends}(\Gamma)=\varprojlim
\pi_0(\Gamma\setminus C)$ is known as the set of \emph{ends} of
$\Gamma$. The number of ends of $\Gamma$ is the cardinality of
$\mathrm{Ends}(\Gamma)$.  The number of ends of a finitely generated
group $G$ is the number of ends of the Cayley-graph of $G$ with respect
to any finite set of generators; this number
is independent of the finite generating set we choose for $G$
and it is either $0$, $1$, $2$ or $\infty$~\cite{St71}.
Here are some examples:
(i) every finite group has $0$ ends, (ii) $\mathbb{Z} \times \mathbb{Z}$
has one end, (iii) $\mathbb{Z}$ has two ends, and (iv) $F_2$ has infinitely
many ends. Stallings' famous ends
theorem~\cite{St71} says that if $G$ is a finitely generated group
with two or more ends, then $G$ splits nontrivially as an
amalgamated product or an HNN-extension over a finite subgroup. This
can be reformulated in terms of actions on trees via Bass-Serre
theory~\cite{Se03}.

A group acts \emph{nontrivially} on a tree if it has no global
fixed-point, i.e., there is no node $v$ in the tree with
$Gv = \{v\}$. A group $G$ is said to \emph{split} over a subgroup $H$ if
there is a nontrivial action of $G$ on a tree $T$ such that $H$ is
the stabilizer of an edge $e$ and the orbit $Ge$ consists of all
edges of $T$.  This is equivalent to $G$ splitting as an amalgamated
product or HNN-extension with $H$ as the amalgamation base,
respectively the associated subgroup~\cite{Se03}.  We shall need the
following simple lemma.

\begin{lemma}\label{l nontrivialaction}
Let $G$ be a finitely generated group with a nontrivial action on a
tree $T$ and let $H\leq G$ be a finite index subgroup.  Then $H$
acts nontrivially on $T$.
\end{lemma}

\begin{proof}
Recall that if $g$ is an automorphism of a tree $T$, then $g$ is
said to be \emph{elliptic} if $g$ fixes some point of $T$.  It is well
known (this follows immediately from~\cite[I.6.4, Proposition~25]{Se03}, 
for instance) that if $g^n$ ($n \geq 1$) is elliptic, then $g$ is elliptic.
Now if $H$ has a global fixed point, then $H$ consists entirely of
elliptic automorphisms of $T$. Let $[G:H]=n$ and $g\in G$.  Then
$g^n\in H$ and hence is elliptic. It follows that every element of
$G$ is elliptic.  But it is well known~\cite[I.6.5, Corollary~3]{Se03}
that any finitely 
generated group of elliptic automorphisms of a tree has a global
fixed point, contradicting that the action of $G$ is nontrivial.  It
follows that the action of $H$ is nontrivial. \qed
\end{proof}

\begin{theorem}\label{t submonoidforinfinitelymanyends}
Let $G$ be a finitely generated group with two or more ends such
that the intersection of all the finite index subgroups of $G$ is
torsion-free. Then $G$ is virtually a free product and hence the
rational subset membership and submonoid membership problems for $G$
are recursively equivalent.
\end{theorem}

\begin{proof}
By Stallings ends theorem~\cite{St71}, $G$ splits nontrivially over
a finite subgroup.  So by Bass-Serre theory~\cite{Se03} $G$ acts
nontrivially on a tree $T$ so that there is one orbit of edges and
the stabilizer of an edge is finite. Let $H$ be an edge stabilizer;
since $H$ is a finite group, by hypothesis there is a normal
subgroup $N\lhd G$ of finite index such $H\cap N=\{1\}$. By
Lemma~\ref{l nontrivialaction} the action of $N$ on $T$ is
nontrivial.  Since each edge stabilizer in $G$ is a conjugate of
$H$, it follows  no element of $N\setminus\{1\}$ fixes an edge. Therefore, $N$
splits nontrivially as a free product~\cite{Se03}. This completes
the proof. \qed\end{proof}

\begin{corollary} \label{C tf or rs}
Let $G$ be a finitely generated group with two or more ends which is
either virtually torsion-free or residually finite.  Then the
rational subset membership and submonoid membership problems for $G$
are recursively equivalent.
\end{corollary}

\begin{proof}
Clearly Theorem~\ref{t submonoidforinfinitelymanyends} applies under
either of these hypotheses. \qed\end{proof}

Let us now come back to graph groups.
Theorems~\ref{T main} and
\ref{T submonoid} imply that the submonoid membership problem
is undecidable for every graph group of the form \[\mathbb{G}(\Sigma \cup \{a\}, I) \simeq
\mathbb{G}(\Sigma, I) \ast \mathbb{Z},\] where
$a \not\in \Sigma$ and $(\Sigma,I)$ is not a transitive forest.
In the rest of the paper, we will sharpen this result. We
show that for a graph group the
submonoid membership problem is decidable if and only if the
rational subset membership problem is decidable, i.e., if and only
if the independence alphabet is a transitive forest. In fact, by our
previous results, it suffices to consider a $\mathsf{P4}$:

\begin{theorem}\label{T P4}
Let $\Sigma= \{a,b,c,d \}$ and $I = \{ (a,b), (b,c), (c,d) \}$, i.e,
$(\Sigma,I)$ is a $\mathsf{P4}$. Then there exists a fixed submonoid
$M$ of $\mathbb{G}(\Sigma,I)$ such that the membership problem of
$M$ within $\mathbb{G}(\Sigma,I)$ is undecidable.
\end{theorem}

\begin{proof}
We follow the strategy of the proof of 
Lemma~\ref{L submonoid free prod}, but working in the graph group
$\mathbb{G}(\Sigma,I)$ makes the encoding more complicated.
Let $R$ denote the trace rewriting system over the trace monoid
$\mathbb{M}(\Sigma^{\pm 1}, I)$ defined in (\ref{system R}),
Section~\ref{SS traces}. As usual denote by $h : (\Sigma^{\pm 1})^* \to
\mathbb{G}(\Sigma,I)$ denote the canonical morphism, which will be
identified with the canonical morphism $h : \mathbb{M}(\Sigma^{\pm
1},I) \to \mathbb{G}(\Sigma,I)$. Let us fix a finite automaton $A$
over the alphabet $\Sigma^{\pm 1}$ such that the membership problem
for $h(L(A))$ within $\mathbb{G}(\Sigma,I)$ is undecidable; such an
automaton exists by Theorem~\ref{T main}. Without loss of generality
assume that 
$$A = (\{ 1, \ldots, n \}, \Sigma^{\pm 1}, \delta, q_0, \{q_f\}),$$
where $\delta \subseteq \{ 1, \ldots, n \} \times (\Sigma^{\pm 1} \cup 
\{ \varepsilon \}) \times \{ 1, \ldots, n \}$
and $q_0 \neq q_f$ (since we allow $\varepsilon$-transitions, we may assume 
that there is only a single final state $q_f$, which is different from the 
initial state $q_0$). 
For a state $q \in  \{ 1, \ldots, n \}$, define the trace
$\widetilde{q} \in \mathbb{M}(\Sigma^{\pm 1}, I)$ by
$$
\widetilde{q} = (ada)^q d (ada)^{-q} = (ada)^q d (a^{-1} d^{-1}a^{-1})^q .
$$
Note that the dependence graph of $\widetilde{q}$ is a linear chain.
Moreover, every symbol from $\Sigma^{\pm 1}$ is dependent on $ad$,
i.e., does not commute with $ad$. The following statement is
straightforward to prove.

\medskip

\noindent {\em Claim 2.} Let $q_1,  \ldots, q_k \in \{1, \ldots,
n\}$, $\varepsilon_1,  \ldots, \varepsilon_k \in \{1, -1\}$ such
that $q_i \neq q_{i+1}$ for all $1 \leq i \leq k-1$. Then
$$
\NF_R(\widetilde{q}_1^{\varepsilon_1} \widetilde{q}_2^{\varepsilon_2} \cdots \widetilde{q}_k^{\varepsilon_k})=
(ada)^{q_1} d^{\varepsilon_1} (ada)^{q_2-q_1}  \cdots
d^{\varepsilon_{k-1}} (ada)^{q_k-q_{k-1}} d^{\varepsilon_k}  (ada)^{-q_k}.
$$
Note that this trace starts (resp. ends) with a
copy of $ada$ (resp. $a^{-1} d^{-1}a^{-1}$).

\medskip

\noindent
Let $\varphi : (\Sigma^{\pm 1})^* \to (\Sigma^{\pm 1})^*$ be the
injective morphism defined by $\varphi(x) = xx$ for $x \in
\Sigma^{\pm 1}$. Thus, $w \in L(A)$ if and only if $\varphi(w) \in
\varphi(L(A))$. Since $(x,y) \in I$ implies that $\varphi(x)$ and
$\varphi(y)$ commute, $\varphi$ can be lifted to an injective morphism $\varphi
: \mathbb{M}(\Sigma^{\pm 1},I) \to \mathbb{M}(\Sigma^{\pm 1},I)$.
The reader can easily verify that, for every trace $t  \in
\mathbb{M}(\Sigma^{\pm 1},I)$, the equality $\NF_R(\varphi(t)) =
\varphi(\NF_R(t))$ holds. In particular, $\varphi(t)$ is irreducible if
and only if $t$ is irreducible and $h(t) = h(u)$ if and only if
$h(\varphi(t)) = h(\varphi(u))$.

Let us fix a trace $t \in \mathbb{M}(\Sigma^{\pm 1},I)$
and define
\begin{equation*}
\Delta = \{ \widetilde{q} \varphi(x) \widetilde{p}^{-1} \mid (q, x, p) \in \delta \} \subseteq
\mathbb{M}(\Sigma^{\pm 1},I) \;\text{ and }\;
u = \widetilde{q}_0 \varphi(t) \widetilde{q}_f^{-1} \in \mathbb{M}(\Sigma^{\pm 1},I).
\end{equation*}
We will show that $h(t) \in h(L(A))$ if and only if $h(u) \in h(\Delta^* )$.

Let us define a
\emph{$1$-cycle} to be a word in $(\Sigma^{\pm 1})^*$of the form
$$
\widetilde{q}_1 \varphi(v_1) \widetilde{q}_2^{-1}\, \widetilde{q}_2 \varphi(v_2) \widetilde{q}_3^{-1}
\cdots \widetilde{q}_{k-1} \varphi(v_{k-1}) \widetilde{q}_k^{-1} \,\widetilde{q}_k 
\varphi(v_k) \widetilde{q}_1^{-1}
$$
such that $k \geq 1$, $q_1, \ldots, q_k \in \{1, \ldots, n\}$,
$v_1, \ldots, v_k \in (\Sigma^{\pm 1})^*$, and
$v_1 \cdots v_k = 1$ in $\mathbb{G}(\Sigma,I)$ (hence, also
$\varphi(v_1) \cdots \varphi(v_k) = 1$ in $\mathbb{G}(\Sigma,I)$). 
Note that a $1$-cycle equals $1$ in $\mathbb{G}(\Sigma,I)$.
We say that a word of the form
$\widetilde{q}_1 \varphi(v_1) \widetilde{p}_1^{-1}\,\widetilde{q}_2 
\varphi(v_2) \widetilde{p}_2^{-1} \cdots
\widetilde{q}_m \varphi(v_m) \widetilde{p}_m^{-1}$,
where $q_1, p_1, \ldots, q_m,p_m \in \{1, \ldots, n\}$ and
$v_1, \ldots, v_m \in (\Sigma^{\pm 1})^*$,
contains a $1$-cycle, if there are positions $1 \leq i \leq j \leq m$
such that $\widetilde{q}_i \varphi(v_i) \widetilde{p}_i^{-1} \cdots
\widetilde{q}_j \varphi(v_j) \widetilde{p}_j^{-1}$ is a $1$-cycle.
If a word does not contain a $1$-cycle, then it is called
\emph{$1$-cycle-free}.

\medskip

\noindent {\em Claim 3.} Let $m \geq 1$ and $$v = \widetilde{q}_1
\varphi(v_1) \widetilde{p}_1^{-1} \; \widetilde{q}_2 \varphi(v_2)
\widetilde{p}_2^{-1} \cdots \widetilde{q}_m \varphi(v_m)
\widetilde{p}_m^{-1},$$ where $q_1, p_1, \ldots, q_m,p_m \in \{1,
\ldots, n\}$ and $v_1, \ldots, v_m \in (\Sigma^{\pm 1})^*$. If $v=1$
in $\mathbb{G}(\Sigma,I)$, then $v$ contains a $1$-cycle.

\medskip

\noindent
{\em Proof of Claim 3.}
We prove Claim 3 by induction over
$m$. Assume that $v=1$ in $\mathbb{G}(\Sigma,I)$. If $m = 1$, then
we obtain the identity
\begin{equation}
\label{L 3}
\widetilde{q}_1 \varphi(v_1) \widetilde{p}_1^{-1} =
(ada)^{q_1} d (ada)^{-q_1} \varphi(v_1)  (ada)^{p_1} d^{-1} (ada)^{-p_1}  = 1
\end{equation}
in $\mathbb{G}(\Sigma,I)$. Assume without loss of generality that
$v_1$, viewed as a trace, is irreducible with respect to $R$. Then
also $\varphi(v_1)$ is irreducible. If $\varphi(v_1) = \varepsilon$
and $p_1 = q_1$, then $v$ is a $1$-cycle. If  $\varphi(v_1) = \varepsilon$,
and $p_1 \neq q_1$, then
we obtain a contradiction, since $\NF_R(\widetilde{q}_1  \widetilde{p}_1^{-1})$
is nonempty by Claim 2. Now assume that
$\varphi(v_1) \neq \varepsilon$.  In the trace
$$
(ada)^{q_1} d (ada)^{-q_1} \varphi(v_1)  (ada)^{p_1} d^{-1} (ada)^{-p_1}
$$
only the last $a^{-1}$ of the factor $(a^{-1} d^{-1} a^{-1})^{q_1}$
may cancel against the first $a$ of $\varphi(v_1)$ (in case $a \in
\min(v_1)$) and the first $a$ of the factor $(ada)^{p_1}$ may cancel
against the last $a^{-1}$ of $\varphi(v_1)$ (in case $a^{-1} \in
\max(v_1)$). To see this, note that if $a \not\in \min(v_1)$, then
$(ada)^{-q_1} \varphi(v_1)$ is irreducible with respect to $R$. If
$a \in \min(v_1)$ then $\varphi(v_1) = aa \varphi(t)$ for some trace
$t$. Then
$$
(a^{-1} d^{-1} a^{-1})^{q_1}  \varphi(v_1) = (a^{-1} d^{-1} a^{-1})^{q_1}
aa\varphi(t) \to_R (a^{-1} d^{-1} a^{-1})^{q_1-1} a^{-1} d^{-1} a \varphi(t).
$$
Since $a$ and $d$ do not commute,
we cannot have $d \in \min(a \varphi(t))$, hence cancellation stops and
$\NF_R((a^{-1} d^{-1} a^{-1})^{q_1}  \varphi(v_1) ) =
(a^{-1} d^{-1} a^{-1})^k a^{-1} d^{-1} a \varphi(t)$ where $k = q_1-1 \geq 0$.
Moreover, if $a^{-1}$ is a maximal symbol of $t$, then $\varphi(t) = \varphi(t')
a^{-1}a^{-1}$ for some trace $t'$. Hence, by making a possible
cancellation with the first $a$ in $(ada)^{p_1}$, it follows finally
that
\begin{equation*}
\NF_R(\widetilde{q}_1 \varphi(v_1) \widetilde{p}_1^{-1}) = (ada)^{q_1} d (ada)^{-k} a^{-1} d^{-1} x d a (ada)^\ell d^{-1} (ada)^{-p_1}
\neq \varepsilon
\end{equation*}
for some trace $x$, where
$\ell = p_1-1 \geq 0$. This contradicts
again (\ref{L 3}) and proves the inductive base case $m = 1$ in
Claim 3.

Now assume that $m  \geq 2$.

\medskip

\noindent
{\em Case 1.} There is $1 \leq i < m$ such that $p_i = q_{i+1}$.
Then $v=1$ in $\mathbb{G}(\Sigma,I)$ implies
\begin{equation*}
\widetilde{q}_1  \varphi(v_1) \widetilde{p}_1^{-1}  \cdots
\widetilde{q}_{i-1}  \varphi(v_{i-1}) \widetilde{p}_{i-1}^{-1}\, \widetilde{q}_i  \varphi(v_i v_{i+1}) \widetilde{p}_{i+1}^{-1}\,
\widetilde{q}_{i+2}  \varphi(v_{i+2}) \widetilde{p}_{i+2}^{-1} \cdots
\widetilde{q}_m  \varphi(v_m) \widetilde{p}_m^{-1}
\end{equation*}
is $1$ in $\mathbb{G}(\Sigma,I)$.
By induction, we can conclude that above word contains
a $1$-cycle. But then also the word $v$ must contain a $1$-cycle.

\medskip

\noindent
{\em Case 2.} $p_i \neq q_{i+1}$ for all $1 \leq i < m$.              
If there is $1 \leq i \leq m$ such that
$v_i = 1$ in $\mathbb{G}(\Sigma,I)$ and $q_i = p_i$ then $v$ contains
the $1$-cycle $\widetilde{q}_i \varphi(v_i) \widetilde{p}_i^{-1}$.
Now assume that $q_i \neq p_i$ whenever $v_i=1$ in  $\mathbb{G}(\Sigma,I)$.
Let $v'$ be the word that results from $v$ by deleting
all factors $\varphi(v_i)$, which are equal $1$ in
$\mathbb{G}(\Sigma,I)$. In the following, we consider $v'$ as a trace.
Consider a maximal factor of $v'$ of the form
\begin{equation}
\widetilde{p}_i^{-1} \widetilde{q}_{i+1}\widetilde{p}_{i+1}^{-1} \widetilde{q}_{i+2} \cdots
\widetilde{p}_{j-1}^{-1} \widetilde{q}_j  \label{abcd-word}
\end{equation}
where  $j \geq i+1$ and $\varphi(v_{i+1}) = \cdots =
\varphi(v_{j-1}) = 1$, $\varphi(v_i)  \neq 1 \neq \varphi(v_j)$
in $\mathbb{G}(\Sigma,I)$. Claim 2 show that the $R$-normal form of this
trace starts (resp. ends) with a copy of
$ada$ (resp. $a^{-1} d^{-1}a^{-1}$), and similarly for maximal prefixes
(resp. suffixes) of the form 
\begin{equation} \label{E pref suf}
\widetilde{q}_1 \widetilde{p}_1^{-1} \cdots
\widetilde{q}_{i-1} \widetilde{p}_{i-1}^{-1} \widetilde{q}_i 
\qquad \text{ (resp. $\widetilde{p}_i^{-1} \widetilde{q}_{i+1}\widetilde{p}_{i+1}^{-1}  \cdots
\widetilde{p}_m^{-1} \widetilde{q}_m$).}
\end{equation}
 In $v'$, factors of the form
(\ref{abcd-word}) and (\ref{E pref suf})
are separated by traces $\varphi(v_i)$, where
$\varphi(v_i) \neq 1$ in $\mathbb{G}(\Sigma,I)$. Without loss of
generality assume that each such trace $\varphi(v_i)$ is irreducible
and hence non-empty. As for the base case
$m=1$, one can show that in such a concatenation, only a single
minimal $a$ and a single maximal $a^{-1}$ of a trace $\varphi(v_i)
\neq \varepsilon$ may be cancelled. It follows that $\NF_R(v) \neq \varepsilon$,
which contradicts $v=1$ in $\mathbb{G}(\Sigma,I)$.
This concludes the proof of Claim 3. 

Now we can prove $h(t) \in h(L(A))$ if and only if $h(u) =
h(\widetilde{q}_0 \varphi(t) \widetilde{q}_f^{-1}) \in h( \Delta^* )$. First
assume that $h(t) \in h(L(A))$. Let $a_1 \cdots a_m \in L(A)$ such
that $(q_{i-1}, a_i, q_i) \in \delta$ for $1 \leq i \leq m$, $q_m =
q_f$, and $a_1 \cdots a_m = t$ in $\mathbb{G}(\Sigma,I)$. Then
\begin{equation*}
h(\widetilde{q}_0 \varphi(t) \widetilde{q}_f^{-1}) = h( \widetilde{q}_0 \varphi(a_1)
\widetilde{q}_1^{-1}  \widetilde{q}_1 \varphi(a_2) \widetilde{q}_2^{-1}
\cdots \widetilde{q}_{m-1} \varphi(a_m) \widetilde{q}_m^{-1} )  \in
h(\Delta^*).
\end{equation*}
Now assume that $h(\widetilde{q}_0 \varphi(t) \widetilde{q}_f^{-1}) \in 
h(\Delta^* )$. Thus,
$$
\widetilde{q}_0 \varphi(t) \widetilde{q}_f^{-1} = \widetilde{q}_1 \varphi(a_1) \widetilde{p}_1^{-1}  \widetilde{q}_2 \varphi(a_2) \widetilde{p}_2^{-1} \cdots
\widetilde{q}_m \varphi(a_m) \widetilde{p}_m^{-1}
$$
in $\mathbb{G}(\Sigma,I)$,
where $q_1, p_1, \ldots, q_m,p_m \in \{1, \ldots, n\}$, $a_1,
\ldots, a_m \in \Sigma^{\pm 1} \cup \{\varepsilon\}$, and $(q_i,
a_i, p_i) \in \delta$ for $1 \leq i \leq m$. Without loss of
generality we may assume that the word 
$\widetilde{q}_1 \varphi(a_1) \widetilde{p}_1^{-1}
\widetilde{q}_2 \varphi(a_2) \widetilde{p}_2^{-1} \cdots \widetilde{q}_m \varphi(a_m) \widetilde{p}_m^{-1}$ is $1$-cycle-free
(otherwise we can remove all $1$-cycles from this word;
note that a $1$-cycle
equals $1$ in the group $\mathbb{G}(\Sigma,I)$).
Let
$$
v = \widetilde{q}_f \varphi(t^{-1}) \widetilde{q}_0^{-1}\,  \widetilde{q}_1 \varphi(a_1) \widetilde{p}_1^{-1}\,  \widetilde{q}_2 \varphi(a_2) 
\widetilde{p}_2^{-1} \cdots \widetilde{q}_m \varphi(a_m) \widetilde{p}_m^{-1}.
$$
Since $v=1$ in $\mathbb{G}(\Sigma,I)$, Claim 3 implies that 
$v$ contains a $1$-cycle. We claim that this $1$-cycle must be the whole word $v$:
first of all, the suffix $\widetilde{q}_1 \varphi(a_1) \widetilde{p}_1^{-1}
\cdots \widetilde{q}_m \varphi(a_m) \widetilde{p}_m^{-1}$ of $v$ is
$1$-cycle-free. If a prefix $\widetilde{q}_f \varphi(t^{-1})
\widetilde{q}_0^{-1}\, \widetilde{q}_1 \varphi(a_1) \widetilde{p}_1^{-1}
\cdots \widetilde{q}_i \varphi(a_i) \widetilde{p}_i^{-1}$ for $i < m$ is a
$1$-cycle, then $\widetilde{q}_{i+1} \varphi(a_{i+1})
\widetilde{p}_{i+1}^{-1} \cdots \widetilde{q}_m \varphi(a_m)
\widetilde{p}_m^{-1} = 1$ in $\mathbb{G}(\Sigma,I)$. Hence, Claim 3 implies that
the word $\widetilde{q}_{i+1} \varphi(a_{i+1}) \widetilde{p}_{i+1}^{-1}
\cdots \widetilde{q}_m \varphi(a_m) \widetilde{p}_m^{-1}$ contains a
$1$-cycle, contradicting the fact that the word $\widetilde{q}_1 \varphi(a_1)
\widetilde{p}_1^{-1} \cdots \widetilde{q}_m \varphi(a_m)
\widetilde{p}_m^{-1}$ is $1$-cycle-free. Thus, indeed, $v$ is a
$1$-cycle. Hence, $q_0 = q_1$, $q_f = p_m$, $p_i = q_{i+1}$ for $1
\leq i < m$, and $t^{-1} a_1 \cdots a_m = 1$ in $\mathbb{G}(\Sigma,I)$, i.e., $h(t) =
h(a_1 \cdots a_m) \in h(L(A))$. This shows that the membership problem
for the submonoid $h(\Delta^*)$ within $\mathbb{G}(\Sigma,I)$ 
is indeed undecidable.     \qed
\end{proof}

Recall that a graph is not a transitive forest if and only if it
either contains an induced $\mathsf{C4}$ or $\mathsf{P4}$
\cite{Wolk65}. Together with Mihailova's result for the generalized
word problem of $F(\{a,b\}) \times F(\{c,d\})$, Theorems~\ref{T main}
and \ref{T P4} imply:

\begin{corollary}\label{C submonoid graph group}
The submonoid membership problem for a graph group
$\mathbb{G}(\Sigma,I)$ is decidable if and only if $(\Sigma,I)$ is a
transitive forest. Moreover, if $(\Sigma,I)$ is not a transitive
forest, then there exists a fixed submonoid $M$ of
$\mathbb{G}(\Sigma,I)$ such that the membership problem for $M$
within $\mathbb{G}(\Sigma,I)$ is undecidable.
\end{corollary}

Since $\mathsf{P4}$ is a chordal graph, the generalized word problem
for $\mathbb{G}(\mathsf{P4})$ is decidable \cite{KaWeMy05}. Hence,
$\mathbb{G}(\mathsf{P4})$ is an example of a group for which the
generalized word problem is decidable but the submonoid membership
problem is undecidable.

\section{Open problems}

The definition of the class $\mathcal C$ at the beginning of 
Section~\ref{S RSMP} leads to the
question whether decidability of the rational subset membership
problem is preserved under direct products with $\mathbb{Z}$.
An affirmative answer would lead in combination with the results
from \cite{KaSiSt06,Ned00} to a more direct proof
of Theorem~\ref{T dec}.

Concerning graph groups, the precise borderline for the decidability
of the generalized word problem remains open. By
\cite{KaWeMy05}, the generalized word problem is decidable if the
independence alphabet is chordal. Since every transitive forest is
chordal, Theorem~\ref{T main} does not add any new decidable cases.
On the other hand, if the independence alphabet contains an induced
$\mathsf{C4}$, then the generalized word problem is
undecidable \cite{Mih66}. But it is open for
instance, whether for a cycle of length 5 the corresponding graph
group has a decidable generalized word problem.

Another open problem concerns the complexity of the rational subset
membership problem for graph groups, where the independence alphabet
is a transitive forest. If the independence alphabet is part of the
input, then our decision procedure does not yield an elementary
algorithm, i.e., an algorithm where the running time  is bounded by
an exponent tower of fixed height. This is due to the fact that each
calculation of the Parikh image of a context-free language leads to
an exponential blow-up in the size of the semilinear sets in the
proof of Lemma~\ref{L free product}. An NP lower bound follows
from the NP-completeness of integer programming.

Theorem~\ref{T submonoid} and \ref{T submonoid HNN}
lead to various research directions. One
might try to get rid of the restriction that $g^{-1} A g \cap A = 1 = 
h^{-1} A h \cap A$ for some $g \in G$, $h \in H$ and the analogous
restrictions for HNN extensions. 
These two results together would imply that Corollary~\ref{C tf or rs}
holds for all groups with two or more ends. 

In fact it is natural to 
ask whether, 
for every finitely generated group $G$, the submonoid membership and
rational subset membership problems are recursively equivalent. By
Theorem~\ref{T submonoid}, this is equivalent to the preservation of
the decidability of the submonoid membership problem under free
products (which is again not known to hold): simply choose for $H$
in Theorem~\ref{T submonoid} any nontrivial group with a decidable
rational subset membership problem. Recall that the decidability of
the generalized word problem as well as the rational subset
membership problem is preserved under free products. Notice that for
a torsion group, the submonoid membership problem is equivalent to the
generalized word problem, while the rational subset membership problem reduces
to membership in products $H_1\cdots H_n$ of finitely generated
subgroups.

%\bibliographystyle{abbrv}
%\bibliography{../../../BIBTEX/bib.bib}

\end{document}